\documentclass[10pt]{amsart}

\usepackage[margin=1.0in]{geometry}
\usepackage{siunitx} 
\RequirePackage[numbers]{natbib}
\RequirePackage[colorlinks,citecolor=blue,urlcolor=blue]{hyperref}

\usepackage{amsthm,amsmath}
\usepackage{amssymb} 
\usepackage{amsfonts,epsfig,epstopdf,url,array}

\usepackage{graphicx}
\graphicspath{{./Figures/}}

\def\acuteangle{{\raisebox{0.05ex}{\includegraphics[height=1.4ex]{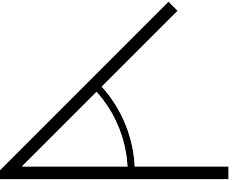}}}}
\def\obtuseangle{{\raisebox{0.05ex}{\includegraphics[height=1.4ex]{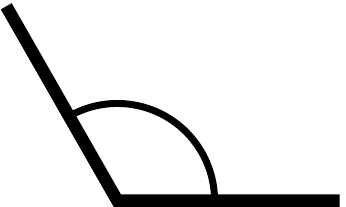}}}}
\def\rightangle{{\raisebox{0.05ex}{\includegraphics[height=1.4ex]{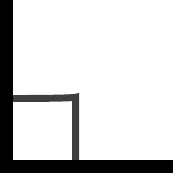}}}}
\def\paragram{{\raisebox{0.05ex}{\includegraphics[height=1.4ex]{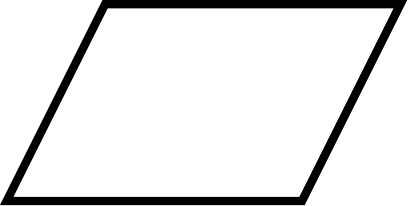}}}}

\usepackage{textcomp}
\usepackage{caption}
\usepackage{longtable}
\newcommand\Tstrut{\makebox[0pt][c]{\rule{0pt}{2.6ex}}}		
\newcommand\Bstrut{\makebox[0pt][c]{\rule[-1.2ex]{0pt}{0pt}}}	

\theoremstyle{remark}

\def\Item$#1${\item {\hfill $\displaystyle#1$ \hfill\refstepcounter{equation}(\theequation)}}

\title{Classifying Rational Parallelepipeds}
\author{Randall L. Rathbun}
\email{randallrathbun@gmail.com}

\begin{document}
\date{November 27, 2018}

\begin{abstract}
It is proposed that the name {\it `Diophantine parallelepiped'} be applied to rational parallelepipeds.
By examining the count of possible right angles {\rightangle} at any vertice, the parallelepipeds become
classified into 5 types: {\it acute} {\it anorthgonal} {\it triclinic}; {\it obtuse} {\it anorthgonal} {\it triclinic};
{\it monorthogonal} {\it bi-clinic}; {\it bi-orthogonal} {\it monoclinic}; and {\it rectangular cuboid}; There is
sufficient confusion in the present literature to warrant an attempt at clarity for these types of pipeds and
to offer a classification scheme.

The possible rational components of the Diophantine parallelepiped are examined and some results from
computer searches are also presented.
\end{abstract}

\maketitle

\section{Introduction to Diophantine Parallelepipeds}

Parallelepipeds have been examined, in Diophantine analysis of Number Theory, but a lack of clarity still
seems to exist.

We define a rational parallelepiped in $_n$-dimensions as a polytope spanned by $n$ vectors
$\boldsymbol{\vec{v}}_1, \dots,\:\boldsymbol{\vec{v}}_n$ in a vector space over the rationals,
$\mathbb{Q}$, or integers, $\mathbb{Z}$
\begin{center}
	span($\boldsymbol{\vec{v}}_1,\dots,\:\boldsymbol{\vec{v}}_n) \: = \: t_1\boldsymbol{\vec{v}}_1 + \: \dots + t_n\boldsymbol{\vec{v}}_n \text{ for } t_i \in \mathbb{Q}$ for $i=1 \to n$
\end{center}
Here we are interested in three dimensions, $n=3$, so a rational parallelepiped is a prism determined by the 3 rational basis vectors
$\boldsymbol{\vec{a}}$, $\boldsymbol{\vec{b}}$, $\boldsymbol{\vec{c}}$. The prism has 8 vertices, 3 pairs of parallel faces
which are all parallelograms, and 12 edges (3 distinct) $\in\mathbb{Q}$.
\begin{equation}
	\text{span}(\boldsymbol{\vec{a}},\boldsymbol{\vec{b}},\boldsymbol{\vec{c}}) \: = \: t_1\boldsymbol{\vec{a}} + t_2\boldsymbol{\vec{b}} + t_3\boldsymbol{\vec{c}} \;\; \text{ for } t_{i=1 \to 3} \in \mathbb{Q}
\label{eq:span}
\end{equation}

\begin{figure}[!h]
\centering
\includegraphics[scale=1.0]{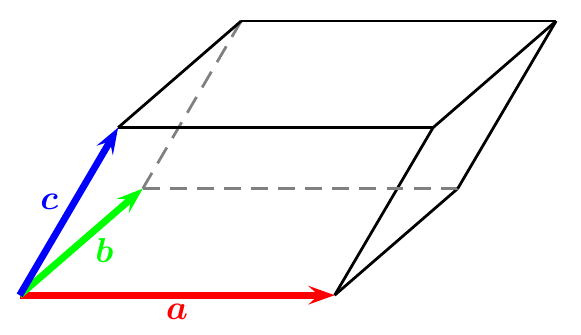}
\captionsetup{justification=centering}
\caption[caption]{The rational parallelepiped with 3 rational basis vectors $\vec{a}$, $\vec{b}$, and $\vec{c}$}
\label{fig:ratpiped}
\end{figure}

In Figure \ref{fig:ratpiped}, the lengths of the basis vectors, $\vec{a}$, $\vec{b}$, and $\vec{c}$, are rational
$\in \mathbb{Q}$, and they determine the parallelepiped.

\section{Classifying the Parallelepipeds by orthogonality}

Literature abounds with examples confusing rectangular parallelepipeds with monoclinic pipeds. The bi-clinic or
mon-orthogonal case is hardly even known. The rectangular parallelepiped is called the `{\it cuboid}', or `{\it brick}',
or `{\it Euler brick}', or even the `{\it orthotope}', engendering confusion over the preferred name for the same object.

In an attempt at clarity for Diophantine analysis, this author proposes a classification scheme based upon the maximum
count of right angles {\rightangle} which might occur at any of the 8 vertices in the Diophantine parallelepiped.

Consider the 3 surface angles from Figure \ref{fig:ratpiped}., $\alpha$ between $\vec{a}$ and $\vec{b}$; $\beta$ between
$\vec{a}$ and $\vec{c}$; and $\gamma$ between $\vec{b}$ and $\vec{c}$. There may be none, or 1, or 2, or 3 {\rightangle} right angles
at a vertice of the Diophantine parallelepiped. All 8 vertices need to be considered.

\section{Five Classes of Rational Diophantine Parallelepipeds}

There are 5 unique classes of parallelepipeds which do exist, depending upon the count
of right angles at the origin [the intersection of the basis vectors], as shown and proved in Appendix A.

\begin{table}[h]
\begin{center}
\caption{The five classes of rational parallelepipeds by surface angles at a vertice.}
\begin{tabular}{c | c | c c c c c c c c c c c}
Name & Comment &
\multicolumn{9}{c}{Vertice Group(s)} \\
\hline
\phantom{1o}acute & (anorthic) [no {\rightangle}] & \texttt{[-1} & \texttt{-1} & \texttt{ 1]} & \texttt{[ }{\color{red}\texttt{1}} & {\color{red}\texttt{ 1}} & {\color{red}\texttt{ 1}}\texttt{]} \\
\phantom{1}obtuse & (anorthic) [no {\rightangle}] & \texttt{[}{\color{red}\texttt{-1}} & {\color{red}\texttt{-1}} & {\color{red}\texttt{-1}}\texttt{]} & \texttt{[-1} & \texttt{ 1} & \texttt{ 1]} \\
1-ortho & (bi-clinic) & \texttt{[-1} & \texttt{-1} & {\color{red}\texttt{ 0}}\texttt{]} & \texttt{[-1} & {\color{red}\texttt{ 0}} & \texttt{ 1]} & \texttt{[ {\color{red}0}} & \texttt{ 1} & \texttt{ 1]} \\
2-ortho & (monoclinic) & \texttt{[-1} & {\color{red}\texttt{ 0}} & {\color{red}\texttt{ 0}}\texttt{]} & \texttt{[ }{\color{red}\texttt{0}} & {\color{red}\texttt{ 0}} & \texttt{ 1]} \\
3-ortho & (rectangular) & \texttt{[ }{\color{red}\texttt{0}} & {\color{red}\texttt{ 0}} & {\color{red}\texttt{ 0}}\texttt{]}
\end{tabular}
\label{table:5classes}
\end{center}
\end{table}

In the table above, the classsifier for the vertice group is shown in the Vertice Group(s) column. The meaning is as follows: $1 =\text{acute } \acuteangle$ angle; $0 = \text{right }\rightangle$ angle; and $-1 =\text{obtuse } \obtuseangle$ angle. The classifier is color highlighted for the right \rightangle { }angle(s).

\noindent
NOTE: The vertice group gives the actual types of vertices which actually occur in the class.

\noindent
This suggests a possible classification scheme with 5 categories or class names:
\begin{itemize}
\item {\it the Diophantine acute anorthic parallelepiped (acute triclinic)}
\item {\it the Diophantine obtuse anorthic parallelepiped (obtuse triclinic)}
\item {\it the Diophantine 1-ortho parallelepiped (biclinic)}
\item {\it the Diophantine 2-ortho parallelepiped (monoclinic)}
\item {\it the Diophantine 3-ortho parallelepiped (rectangular)}
\end{itemize}
Please note that for the first two classes, the {\it acute anorthic}, and the {\it obtuse anorthic},
that no right angle {\rightangle} exists at any vertex of the parallelepiped.

The labeling of the pipeds using the terms anorthic, or triclinic, the relatively unknown biclinic and the well known
monoclinic is deliberate, and comes from the mineral crystallography classification scheme, since the crystals systems share
some of the same morphology as the Diophantine pipeds.

Thus the proposed classification scheme then depends upon whether or not (0,1,2,3) right angles exist for at least one
of the eight vertices of the piped (actually at the origin or intersection of the basis vectors creating the piped).

\begin{table}
\begin{center}
\caption{Diophantine parallelepiped classification scheme into 5 classes.}
\begin{tabular}{|c|c|c|c|c|c|}
\hline
{\small proposed name} & {\small\it acute} & {\small\it obtuse} & {\small\it 1-ortho} & {\small\it 2-ortho} & {\small\it rectangular} \\
\hline
{\small common name} & {\small triclinic} & {\small triclinic} & {\small parallelepiped} & {\small monoclinic} & {\small cuboid} \\
\hline
{\small right angle(s) \rightangle } & 0 & 0 & 1 & 2 & 3 \\
\hline
\end{tabular}
\end{center}
\label{table:5scheme}
\end{table}

\section{Rational Parallelepipeds and their components}

We need to take a careful look at the rational components of the Diophantine piped for Diophantine analysis.
Due to the extensive symmetry and parallelism of the parallelepiped figure, many of its components
are congruent, thus insuring rationality for all elements in the same congruent group, such as the edges,
the face diagonals, the skew triangles, the face or body triangles, and the parallelogram areas.

\subsection{Derivation from 3 rational or integer basis-vectors}

\begin{figure}[!ht]
\centering
\includegraphics[scale=1.0]{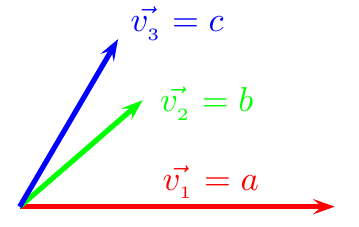}
\caption{The three basis vectors, $\vec{v_1}$, $\vec{v_2}$ $\vec{v_3}$.\\with magnitudes $a$, $b$, $c$ $\in\mathbb{Z}$}
\label{fig:v1v2v3}
\end{figure}

Again, as in eq(\ref{eq:span}) we first start with 3 finite vectors which span $\mathbb{R}^3$ as shown
in figure \ref{fig:v1v2v3}. Let their magnitudes be rational $\in\mathbb{Q}$ or preferably as integer $\in\boldsymbol(Z)$.

Call the integer lengths of these vectors $a$, $b$, and $c$ respectively. Using these 3 vectors
we can create 3 unique parallelograms, using a pair of the 3 vectors, and then combine
the 3 parallelograms into a parallelepiped.

Where 3 parallelograms exist, with their sides of length $[a,b]$, and $[a,c]$, and $[b,c]$, then a parallelepiped
solution exists, where $a$, $b$, and $c$ are the magnitudes of the 3 basis vectors. We label the 3 pairs of the diagonals
of the 3 parallelograms as {\it d,e}; {\it f,g}; and {\it h,j}.
\begin{table}[!ht]
\centering
\begin{tabular}{c | c | c}
parallelogram & sides & diagonals \\
\hline
1 & {\it a,b} & {\it d,e} \\
2 & {\it a,c} & {\it f,g} \\
3 & {\it b,c} & {\it h,j}
\end{tabular}
\caption{3-basis vector solution $a,b,c$ \\ for 3 parallelograms with matching sides}
\label{table:3vect}
\end{table}

\subsection{The Parallelogram equation and rational parallelograms}

The parallelogram equation (\ref{eq:paraeq}) shows that we can pick 3 values of a parallelogram to be rational $\in\mathbb{Q}$.
\begin{equation}
2 (a^2 + b^2) = c^2 + d^2
\label{eq:paraeq}
\end{equation}
where $a,b$ are the sides and $c,d$ are the diagonals of the parallelogram. We can pick 3 rational values for a parallelogram, satisfying
equation (\ref{eq:paraeq}), say sides $a,b$ and the diagonal $c$, then the other diagonal $d$ may or may not be rational. If the second diagonal is rational, we call the parallelogram the {\it rational} or {\it integer} {\it parallelogram} (upon scaling).

It is noted that these integer parallelograms are created from integer triangles, so every parallelogram has at least 1 rational diagonal, and
if a parallelepiped is constructed from 3 such triangles, then at least 3 diagonals, say $d$, $f$, and $h$, from Table \ref{table:3vect} are rational,
while the remaining three $e$, $g$, and $j$ may or may not be rational.

\subsection{Locating coordinates in $\mathbb{R}^3$ for permuted tetrahedrons of a piped}

Using these 3 parallelograms, which have at least 1 rational diagonal each, 48 tetrahedrons can be created from the permutations,
as 6 sets of 8 tetrahedrons. However only 1 set of the 8 tetrahedrons is necessary to completely specify all 6 sets.

A convenient and unique collection covering the 8 possible arrangements of a convenient set is shown in Table \ref{table:tetras} below.

\begin{table}[!h]
\centering
\caption{Combinations of the sides/diagonals of 3 Parallelograms for a Tetrahedron}
\begin{tabular}{c|c|c|c}
\multicolumn{4}{c}{Tetrahedron Sides} \\
\hline
{\texttt a d b f c h} & {\texttt a d b f c j} & {\texttt a d b g c h} & {\texttt a d b g c j} \\
{\texttt a e b f c h} & {\texttt a e b f c j} & {\texttt a e b g c h} & {\texttt a e b g c j}
\end{tabular}
\label{table:tetras}
\end{table}

Depending upon whether or not the second diagonal(s) is rational or not, we can create
an {\it integer tetrahedron}, in at least 1 case of the 8 given above.

The reason we deal with {\it integer tetrahedrons} is that they provide us with the 6 key lengths: $a$, $b$, $c$, $d$, $e$, $f$; which
enable us to fix the vertices of the Parallelepiped into $\mathbb{R}^3$ Cartesian space.

An {\it integer tetrahedron} is shown below in figure \ref{fig:tetra} with 6 integer (or rational) sides,
$a$, $b$, $c$, $d$, $e$, \& $f$. These sides are slightly relabeled from the set of 8 previously given so
that the equations (\ref{eq:4}, \ref{eq:5}) below giving the lattice locations are correct.

You can quickly recognize the 3 vectors $\vec{v_1}$, $\vec{v_2}$ and $\vec{v_3}$ creating edges $a$, $c$, and $e$
in the figure:

\begin{figure}[!ht]
\centering
\includegraphics[scale=1.0]{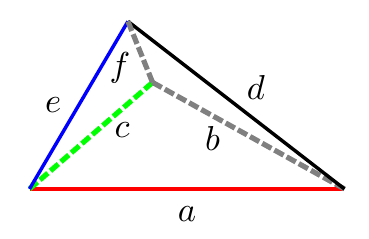}
\caption{Integer tetrahedron.}
\label{fig:tetra}
\end{figure}

We have to locate the vertices of this tetrahedron, given the lengths, into $\mathbb{R}^3$, and I choose
the following method to fix the vertices $v_1$, $v_2$, $v_3$, and $v_4$ into Cartesian 3D space. We find the
3 vertices first.
\begin{equation}
\begin{aligned}
\text{let tetrahedron edges be } & \; a,b,c,d,e,f \in Z \\
v_1 & = [0,0,0] \\
v_2 & = [a,0,0] \\
\text{let }	\cos(B) & = \frac{a^2+c^2-b^2}{2ac} \\
\text{and }	\sin(B) & = \frac{\sqrt{(a+b+c)(a+b-c)(a+c-b)(b+c-a)}}{2ac} \\
v_3 & = c[\cos(B), \sin(B),0] \\
\end{aligned}\label{eq:4}
\end{equation}
By using the intersections of three $S^2$ spheres from the 3 known vertices $[0,0,0]$, $v_1$, and $v_2$ with the fourth vertice $v_4$, that vertice can be located.
\begin{equation}
\begin{aligned}
	x & = \frac{a^2+e^2-d^2}{2a} \\
	y & = \frac{c^2 + e^2 - f^2 - 2rx}{2s} \\
\text{where } r & = c\cdot\cos(B) \\
\text{and }	s & = c\cdot\sin(B) \text{ both from } v_3 \\
\text{then } z & = \pm \sqrt{e^2 - x^2 - y^2} \\
	v_4 & = [x,y,z]
\end{aligned}\label{eq:5}
\end{equation}

The $+z$ value was used (although $-z$ is another solution).

Most of the time, the vertices are located in a quadratic field $K = \boldsymbol{Q}(\sqrt{D})$, occasionally they can be on
a rational lattice.

We have to follow the same process in equations (\ref{eq:4}),(\ref{eq:5}) for all 8 {\it integer tetrahedrons}
in Table \ref{table:tetras}.
The parallelepiped (1 of 8) is derived from these 3 vectors and we strictly label the vertices as shown and follow this order
for the rest of this paper.

\begin{figure}[!ht]
\centering
\includegraphics[scale=1.0]{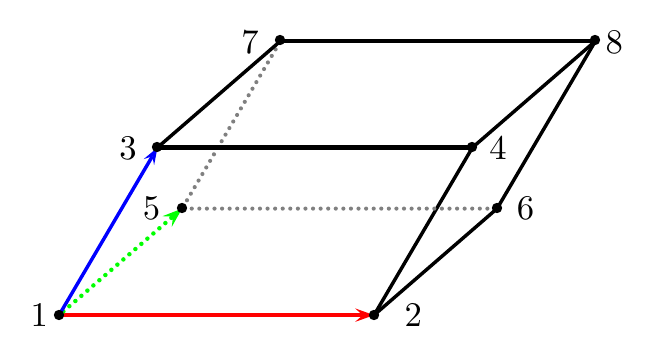}
\caption{The labeled Diophantine parallelepiped.}
\label{fig:vertex}
\end{figure}

This labeled parallelepiped contains several (possibly rational) subcomponents, and they are examined in
the following \S\S.

\subsection{The 28 Diophantine Piped Diagonals}

The first subcomponent to be considered are the diagonals, a line segment between 2 vertices.
There are 28 diagonals in the Diophantine piped, the parallelogram in which they reside is
given also, if it exists.

\begin{center}
\begin{longtable}{c c c c | c c }
\caption{28 Diagonals in the Diophantine piped with parallelograms} \\
\multicolumn{4}{c}{12 Face Diagonals} & \paragram & \paragram \\
\hline
1-4 & 2-3 & 5-8 & 6-7 & F1234 & F5678 \\
1-6 & 2-5 & 3-8 & 4-7 & F1256 & F3478 \\
1-7 & 3-5 & 2-8 & 4-6 & F1357 & F2468 \\
\\[-6pt]
\multicolumn{4}{c}{4 Body Diagonals} & \paragram & \paragram \\
\hline
1-8 & 3-6 & 2-7 & 4-5 & B1368 & B2457 \\
 & & & & B1278 & B3456 \\
 & & & & B1458 & B2367 \\
\newpage
\multicolumn{4}{c}{12 edges} & \paragram & \paragram \\
\hline
1-2 & 3-4 & 5-6 & 7-8 & F1234 & F5678 \\
1-3 & 2-4 & 5-7 & 6-8 & F1234 & F5678 \\
1-5 & 2-6 & 3-7 & 4-8 & F1256 & F3478
\label{table:28diags}
\end{longtable}
\end{center}

There are 6 face parallelograms containing 12 diagonals,
and 6 body parallelograms with 4 unique diagonals, since they share the same set of 4
diagonals in a permutated manner.

\subsection{The 56 Diophantine Piped Triangles}

Triangles are composed of 3 connected diagonals, or any 3 vertices.
There are 56 triangles to consider in the Diophantine piped. However only 48 reside
on the face or body parallelograms.

\begin{center}
\begin{longtable}{c c c | c c c}
\caption{48 Diophantine Piped Triangles on Parallelograms.} \\
$\bigtriangleup$ & $\bigtriangleup$ & \paragram & $\bigtriangleup$ & $\bigtriangleup$ & \paragram \\
\hline
123 & 234 & F1234 & 567 & 678 & F5678 \\
124 & 134 & F1234 & 568 & 578 & F5678 \\ [0.5em]
125 & 256 & F1256 & 347 & 478 & F3478 \\
126 & 156 & F1256 & 378 & 348 & F3478 \\ [0.5em]
135 & 357 & F1357 & 246 & 468 & F2468 \\
137 & 157 & F1357 & 248 & 268 & F2468 \\ [0.5em]
127 & 278 & B1278 & 345 & 456 & B3456 \\
128 & 178 & B1278 & 346 & 356 & B3456 \\ [0.5em]
136 & 368 & B1368 & 245 & 457 & B2457 \\
138 & 168 & B1368 & 247 & 257 & B2457 \\ [0.5em]
148 & 158 & B1458 & 236 & 367 & B2367 \\
145 & 458 & B1458 & 237 & 267 & B2367
\label{table:ptris}
\end{longtable}
\end{center}

There are also 8 triangles unaccounted for. These are not found
on the parallelogram components of the Diophantine piped. I call
them the {\it `skew triangles'} as they are not located on any of the
12 parallelograms.

The 8 skew non-parallelogram triangles are

\begin{center}
\begin{longtable}{c c c c}
\caption{8 Skew Triangles} \\
$\bigtriangleup$ & $\bigtriangleup$ & $\bigtriangleup$ & $\bigtriangleup$ \\
\hline
146 & 147 & 167 & 467 \\
853 & 852 & 832 & 532
\label{table:skew}
\end{longtable}
\end{center}

The matching congruent triangles listed above are given in a column,
preserving the isometry of the lengths of any two vertices of the triangle,
following the order of the vertice labeling.

\subsection{The 12 Diophantine Piped Parallelograms}

The Diophantine piped also contains 12 parallelograms.

\begin{center}
\begin{longtable}{c c c c c c}
\caption{12 Parallelograms in the Piped.} \\
\multicolumn{6}{c}{face parallelograms} \\
\multicolumn{2}{c}{\paragram} & \multicolumn{2}{c}{\paragram} & \multicolumn{2}{c}{\paragram} \\
\hline
\multicolumn{2}{c}{F1234} & \multicolumn{2}{c}{F1357} & \multicolumn{2}{c}{F1256} \Tstrut \\
\multicolumn{2}{c}{F5678} & \multicolumn{2}{c}{F2468} & \multicolumn{2}{c}{F3478} \\
\\[-9pt]
\multicolumn{6}{c}{body parallelograms} \\
\paragram & \paragram & \paragram & \paragram & \paragram & \paragram \\
\hline
B1278 & B1368 & B1458 & B2367 & B2457 & B3456 \Tstrut
\label{table:12para}
\end{longtable}
\end{center}

The face parallelograms given above in a column in Table \ref{table:12para}. are congruent, and have
the same area and diagonal lengths. However please note that the 6 body parallelograms have distinct areas.

\subsection{The Diophantine Parallelepiped Itself}

Finally there is the parallelepiped itself, which has 8 ordered
vertices such that there is a pair of 4 vertices which are co-planar
and parallel. This figure is conventionally called the prism.

There is one component to consider for rationality, the volume: \par
\begin{table}[!ht]
\centering
\begin{tabular}{c}
parallelepiped volume \\
\hline
P12345678 \Tstrut
\end{tabular}
\end{table}

\subsection{Summarizing the Diophantine Analysis for all 83 Components}

Considering all the components and the duplications of lengths which occur in the piped, there is a unique
signature which can be checked for rationality:

\noindent
{\bf Rational Signature Considerations}
\begin{itemize}
\item  3 edge lengths (already rational)
\item  6 face diagonals (3 pairs matches lengths)
\item  4 body diagonals
\item 56 triangle areas* (48 subsumed in parallelograms, 8 are not)
\item  4 skew triangle areas
\item  3 face parallelograms diagonals and area
\item  6 body parallelograms diagonals and area
\item  1 piped volume
\end{itemize}

We don't actually check all 56 triangles*, but let the 48 triangles on the parallelograms
be subsumed by the rationality checks upon the body or face parallelograms. However, we must check
4 of the 8 skew triangles for rational area.

Additionally, the volume definitely needs to be checked, as it was discovered by computer testing
that certain combinations of $a$, $b$, $c$, $d$, $e$, and $f$ for the integer tetrahedron
edges actually led to degenerate `flat' tetrahedrons where all three vectors
$\vec{v_1}$, $\vec{v_2}$, $\vec{v_3}$ were co-planer, or spanning only $\mathbb{R}^2$.

Putting all this together results in the unique rational check signature for a given Diophantine parallelepiped.

\subsection{Unique Rationality Check of the Diophantine Piped}

We discover that not all 83 components need to be checked for rationality. Since there are congruences between the lengths
of the edges and diagonals, and congruences between the areas of the triangles, face parallelograms, and body parallelograms
then not all need to be checked. In fact:
\begin{itemize}
\item {\it 27 components need to be checked for rationality}
\end{itemize}
\noindent
We provide an example of this signature for three sample pipeds, where 0 denotes an irrational value and 1 a rational value.
\begin{table}[!ht]
\centering
\caption{Three sample Diophantine parallelepipeds}
{\scriptsize
\begin{tabular}{c c c | c c c c c c}
\multicolumn{3}{c}{Basis Vectors} & \multicolumn{6}{c}{Tetrahedron Sides} \\
$v_1$ & $v_2$ & $v_3$ & a & b & c & d & e & f \\
\hline
\\[-0.7em]
(103,0,0) & ($\frac{-3179}{103}$,$\frac{120}{103}\sqrt{17458}$,0) & ($\frac{-3179}{103}$,$\frac{-2464734}{4495435}\sqrt{17458},\frac{1122}{43645}\sqrt{27915342}$) & 103 & 204 & 157 & 204 & 157 & 264 \Tstrut \Bstrut \\
(10,0,0) & (-5,$6\sqrt{14}$,0) & (0,$\frac{-24}{7}\sqrt{14}$,$\frac{24}{7}\sqrt{35}$) & 10 & 27 & 23 & 26 & 24 & 41 \Tstrut \Bstrut \\
(44,0,0) & (0,117,0) & (0,0,240) & 44 & 125 & 117 & 244 & 240 & 267 \Tstrut
\end{tabular}
}
\label{table:sigs}
\end{table}

\begin{table}[!ht]
\centering
\caption{27 rational component checks for the three sample Diophantine pipeds}
{\footnotesize
\begin{tabular}{| c c c | c c c c | c c c c c c | c c c c | c c c | c c c c c c | c |}
\multicolumn{3}{c}{edges} &
\multicolumn{4}{c}{skew triangles} &
\multicolumn{6}{c}{face diagonals} &
\multicolumn{4}{c}{body diagonals} &
\multicolumn{3}{c}{face area} &
\multicolumn{6}{c}{body area} &
\multicolumn{1}{c}{vol} \\
1 & 1 & 1 & 1 & 0 & 1 & 0 & 1 & 1 & 1 & 1 & 1 & 1 & 0 & 0 & 0 & 0 & 0 & 0 & 1 & 0 & 0 & 0 & 0 & 0 & 0 & 0 \\
1 & 1 & 1 & 0 & 0 & 0 & 0 & 1 & 1 & 1 & 1 & 1 & 1 & 1 & 1 & 1 & 0 & 1 & 0 & 0 & 0 & 0 & 0 & 0 & 0 & 0 & 0 \\
1 & 1 & 1 & 0 & 0 & 0 & 0 & 1 & 1 & 1 & 1 & 1 & 1 & 0 & 0 & 0 & 0 & 1 & 1 & 1 & 1 & 1 & 1 & 1 & 1 & 1 & 1
\end{tabular}
}
\label{table:signature}
\end{table}

\section{Rational constraints for Diophantine objects}

There is an geometric increase in the rationality requirements as the embedding dimension of the geometric object increases.
For this reason, it is necessary to examine both rational triangles and parallelograms, before considering
{\it Diophantine parallelepipeds}, in order to better understand those requirements.

\subsection{Rational Triangles}

In Diophantine analysis of Number Theory, we can study triangles, for example, for rational components. Since the triangle
is 3 sided polygon, (it is actually determined by 2 non-linear vectors), the following properties need to be checked
for rationality.

\noindent
\textbf{Triangle Constraints}
\begin{itemize}
\item rational sides (3 sides)
\item rational area (1 area)
\end{itemize}

This means that $2^4 = 16$ possibilities exist, when considering the rationality of the triangle. But collectively all 3 sides are
considered rational and we definitely are interested if the area is rational. In fact triangles with this property have
been given a special name: \textit{`Heron triangles'}.

We won't consider internal parts of the triangle, such as the medians, or the angle bisectors, or the perpendicular bisectors,
or other civians, for example, although their study in Diophantine analysis is important. We are simply making the point
here that we have 2 components to consider for rationality: the {\it sides} and the {\it area}.

If we force the triangle to have 3 integer sides, $a,b,c \in \mathbb{Q}$, then we only have to consider the area $\in\mathbb{Q}$.

The answer to the question, {\it can we find a triangle with integer sides and rational area?} has been known in antiquity,
the 3,4,5 right triangle with area 6 is the classic answer.

Less well known is the fact that the primitive 5,5,6 isosceles Heron triangle, and the primitive 5,5,8 isosceles Heron triangle,
both with area 12, are the next smallest integer triangles with rational area. In figure \ref{fig:RT-H-H} below, note that the
5-5-6 Heron triangle was created from joining two 3-4-5 right triangles along the side length 4.

\begin{figure}[!ht]
\centering
\includegraphics[scale=1.0]{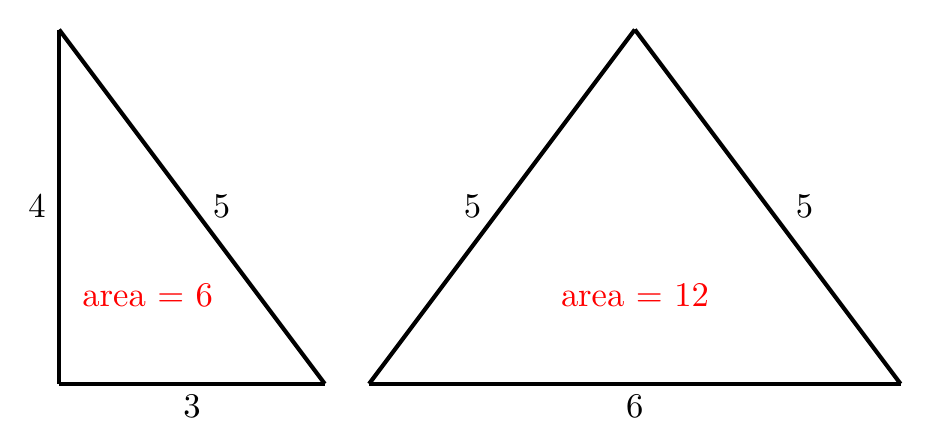}
\captionsetup{justification=centering}
\caption[caption]{The 3,4,5 Right Triangle and 5,5,6 isosceles Heron Triangle}
\label{fig:RT-H-H}
\end{figure}

The next 2 triangles are the composite 6,8,10 right triangle with area 24, and the primitive Heron triangle 9,10,17 with
area 36 which is not isosceles. It should be pointed out that if the triangle sides $\in \mathbb{Z}$, and the triangle
has rational area, then its area $\in \mathbb{Z}$ also. If the sides are fractional $\in \mathbb{Q}$, and if the area
is rational, then the area $\in \mathbb{Q}$.

\noindent Of note is the Heron formula for area of a triangles with sides $a$, $b$, and $c$
\begin{equation}
A = \frac{1}{4}\sqrt{(a+b+c)(-a+b+c)(a-b+c)(a+b-c)}
\end{equation}
which is used to help discover rational area for the triangle in question. These triangles can be parameterized, see \S \ref{sub:params}.

\subsection{Rational Parallelograms}\label{sub:ratpara}

We now consider parallelograms. They are 4 sided planar figures, with 2 pairs of congruent parallel edges;
($a$, $c$), $a=c \text{ and } a \parallel c $ and ($b$, $d$), $b = d \text{ and } b \parallel d$. They also can be
determined by two non-linear vectors but there are 2 additional considerations to be added to the Diophantine analysis.

\noindent The parallelogram has the following properties to check for rationality.

\noindent
\textbf{Parallelogram Constraints}
\begin{itemize}
\item rational sides (4 sides)
\item rational diagonals (2 diagonals)
\item rational area (1 area)
\end{itemize}

This means that we have to consider any one of the $2^7 = 128$ permutations of the parallelograms properties for rationality.
Let the parallelogram have sides ($a$,$c$), ($b$,$d$), with diagonals $d_1$,$d_2$ and area $A$. Considering the pairs of
sides as a unit, we then have only 32 ($2^5$) combinations to consider. But forcing both pairs of sides to be rational further
shortens down the 128 possibilities to just 8 cases.

\begin{table}[!ht]
\centering
\caption{Possible rationality checks for a parallelogram}
\begin{tabular}{ r c | c c | c c | l}
parallelogram & component & $(a,c)$ & $(b,d)$ & $d_1$ & $d_2$ & $A$ \\
\hline
case 1 & rational ? & \checkmark & \checkmark & \textendash & \textendash & \textendash \\
case 2 & rational ? & \checkmark & \checkmark & \checkmark & \textendash & \textendash \\
case 3 & rational ? & \checkmark & \checkmark & \textendash & \checkmark & \textendash \\
case 4 & rational ? & \checkmark & \checkmark & \checkmark & \checkmark & \textendash \\
case 5 & rational ? & \checkmark & \checkmark & \textendash & \textendash & \checkmark \\
case 6 & rational ? & \checkmark & \checkmark & \checkmark & \textendash & \checkmark \\
case 7 & rational ? & \checkmark & \checkmark & \textendash & \checkmark & \checkmark \\
case 8 & rational ? & \checkmark & \checkmark & \checkmark & \checkmark & \checkmark \\
\hline
\end{tabular}
\end{table}

From examination, it can be seen that case 3 is identical to case 2 and case 7 is identical
to case 6, because order is not important, so this narrows our possibilities to just 6 actual cases:

\begin{table}[!ht]
\centering
\caption{Actual rationality checks for the 6 cases of a parallelogram}
\begin{tabular}{ r c | c c | c c | l}
parallelogram & component & $(a,c)$ & $(b,d)$ & $d_1$ & $d_2$ & $A$ \\
\hline
case 1 & rational ? & \checkmark & \checkmark & \textendash & \textendash & \textendash \\
case 2 & rational ? & \checkmark & \checkmark & \checkmark & \textendash & \textendash \\
case 3 & rational ? & \checkmark & \checkmark & \checkmark & \checkmark & \textendash \\
case 4 & rational ? & \checkmark & \checkmark & \textendash & \textendash & \checkmark \\
case 5 & rational ? & \checkmark & \checkmark & \checkmark & \textendash & \checkmark \\
case 6 & rational ? & \checkmark & \checkmark & \checkmark & \checkmark & \checkmark \\
\hline
\label{table:para-rat-checks}
\end{tabular}
\end{table}

Mathematicians would be interested in case 2 and definitely in case 3. In fact case 3 has been
recently parameterized\cite{wyss1,wyss2} (see \S \ref{sub:params}). They might consider case 4
an annoyance, to be avoided. Case 5 might be momentarily considered, but then abandoned because it
has 1 irrational diagonal. So this leaves case 6 as the one which would definitely be considered,
since it has all components rational.

I have purposefully gone through this exercise for both triangles and parallelograms, to emphasize
that the possible rational cases have to be carefully delineated, when studying the Diophantine object in question.
Sorry to say, sometimes this careful analysis is often overlooked, particularly when considering rational parallelepipeds.

As you can imagine, the rational components increase geometrically, when the embedding dimension containing
the object increments, because the figure becomes more complex, and it become important to account for all the
possibilities.

\subsection{Some examples of the smallest integer parallelograms}

The following parallelograms are the smallest integer parallelograms that exist (the next 2 tables list the parallelograms
two to a row).

\begin{table}[!ht]
\centering
\caption{Case 3. First 10 parallelograms with 2 rational diagonals.}
\begin{tabular}{r r | r r || r r | r r}
\multicolumn{2}{c}{sides} & \multicolumn{2}{c}{diags} & \multicolumn{2}{c}{sides} & \multicolumn{2}{c}{diags} \\
\hline
3 & 4 & 5 & 5 & 4 & 7 & 7 & 9 \\
5 & 5 & 6 & 8 & 5 & 10 & 9 & 13 \\
5 & 12 & 13 & 13 & 6 & 7 & 7 & 11 \\
6 & 8 & 10 & 10 & 6 & 13 & 11 & 17 \\
6 & 17 & 17 & 19 & 7 & 9 & 8 & 14
\end{tabular}
\label{table:10-2-diags}
\end{table}

\begin{table}[!ht]
\centering
\caption{Case 6. First 10 parallelograms with 2 diagonals rational and area rational.}
\begin{tabular}{r r | r r || r r | r r}
\multicolumn{2}{c}{sides} & \multicolumn{2}{c}{diags} & \multicolumn{2}{c}{sides} & \multicolumn{2}{c}{diags} \\
\hline
3 & 4 & 5 & 5 & 5 & 5 & 6 & 8 \\
5 & 12 & 13 & 13 & 6 & 8 & 10 & 10 \\
7 & 24 & 25 & 25 & 8 & 15 & 17 & 17 \\
9 & 12 & 15 & 15 & 9 & 40 & 41 & 41 \\
10 & 10 & 12 & 16 & 10 & 24 & 26 & 26
\end{tabular}
\label{table:10-heron}
\end{table}
Please note that if the 2 diagonals are equal, the parallelogram is a rectangle or square. There are some
parallelograms in common with both lists.

\newpage

\begin{table}[!ht]
\centering
\caption{Parallelograms in common - 4 right triangles - 1 Heron triangle.}
\begin{tabular}{r r | r r || r r | r r}
\multicolumn{2}{c}{sides} & \multicolumn{2}{c}{diags} & \multicolumn{2}{c}{sides} & \multicolumn{2}{c}{diags} \\
\hline
3 & 4 & 5 & 5 & 5 & 12 & 13 & 13 \\
5 & 5 & 6 & 8 & 6 & 8 & 10 & 10 \\
7 & 24 & 25 & 25 & & & &
\end{tabular}
\label{table:in_common}
\end{table}

\subsection{Some examples of rational tetrahedrons}\label{sub:rattetra}

From these smallest integer parallelograms, come the smallest integer tetrahedrons. A computer search quickly
discovers the following integer tetrahedrons.
\begin{table}[!ht]
\centering
\caption{The first 40 smallest integer tetrahedrons.}
\begin{tabular}{r r r r r r | r r r r r r}
\multicolumn{6}{c}{tetrahedron edges} & \multicolumn{6}{c}{tetrahedrons edges} \\
\hline
5 & 6 & 5 & 9 & 10 & 13 & 5 & 6 & 5 & 9 & 10 & 9 \\
5 & 6 & 5 & 13 & 10 & 13 & 5 & 6 & 5 & 13 & 10 & 9 \\
5 & 6 & 5 & 13 & 12 & 13 & 5 & 8 & 5 & 9 & 10 & 13 \\
5 & 8 & 5 & 9 & 10 & 9 & 5 & 8 & 5 & 13 & 10 & 13 \\
5 & 8 & 5 & 13 & 10 & 9 & 5 & 8 & 5 & 13 & 12 & 13 \\
7 & 8 & 9 & 21 & 22 & 17 & 7 & 8 & 9 & 21 & 22 & 29 \\
7 & 8 & 9 & 25 & 22 & 17 & 7 & 8 & 9 & 25 & 22 & 29 \\
7 & 12 & 11 & 13 & 16 & 15 & 7 & 12 & 11 & 13 & 16 & 23 \\
7 & 12 & 11 & 21 & 16 & 15 & 7 & 12 & 11 & 21 & 16 & 23 \\
7 & 13 & 16 & 21 & 22 & 18 & 7 & 13 & 16 & 21 & 22 & 34 \\
7 & 13 & 16 & 25 & 22 & 18 & 7 & 13 & 16 & 25 & 22 & 34 \\
7 & 14 & 9 & 21 & 22 & 17 & 7 & 14 & 9 & 21 & 22 & 29 \\
7 & 14 & 9 & 25 & 22 & 17 & 7 & 14 & 9 & 25 & 22 & 29 \\
7 & 14 & 11 & 13 & 16 & 15 & 7 & 14 & 11 & 13 & 16 & 23 \\
7 & 14 & 11 & 21 & 16 & 15 & 7 & 14 & 11 & 21 & 16 & 23 \\
7 & 21 & 16 & 21 & 22 & 18 & 7 & 21 & 16 & 21 & 22 & 34 \\
7 & 21 & 16 & 25 & 22 & 18 & 7 & 21 & 16 & 25 & 22 & 34 \\
7 & 21 & 22 & 25 & 24 & 26 & 7 & 21 & 22 & 25 & 24 & 38 \\
7 & 25 & 22 & 25 & 24 & 26 & 7 & 25 & 22 & 25 & 24 & 38 \\
13 & 10 & 13 & 85 & 84 & 85 & 13 & 24 & 13 & 85 & 84 & 85
\end{tabular}
\label{table:smallest-tetra}
\end{table}

\subsection{The Tetrahedron Permutation Group}\label{sub:tetraperm}

It is important to note that for integer tetrahedrons, using just the six side lengths $a$, $b$, $c$, $d$, $e$, \& $f$,
that up to 24 tetrahedrons can be found in a family. The family is found by considering all permutations of sides,
and accounting for those which create a valid tetrahedron. This permutation group accounts for the appearance
of the similar tetrahedrons seen in Table \ref{table:smallest-tetra} above.

\begin{center}
\begin{longtable}{| c c c c c c |}
\caption{The Tetredron Permutation Group} \\
\hline
\endfirsthead
\multicolumn{6}{c}
{\tablename\ \thetable\ -- The Tetradron Permutation Group \textit{ - continued}} \\
\hline
\endhead
\hline \multicolumn{6}{r}{\textit{continued on next page}} \\
\endfoot
\hline
\multicolumn{6}{l}{The actual group is a set union of the 24 rows given above.} \\
\endlastfoot
$a$ & $b$ & $c$ & $d$ & $e$ & $f$ \\
$a$ & $b$ & $c$ & $\sqrt{2(a^2+e^2)-d^2}$ & $e$ & $\sqrt{2(c^2+e^2)-f^2}$ \\
$a$ & $\sqrt{2(a^2+c^2)-b^2}$ & $c$ & $d$ & $e$ & $\sqrt{2(c^2+e^2)-f^2}$ \\
$a$ & $\sqrt{2(a^2+c^2)-b^2}$ & $c$ & $\sqrt{2(a^2+e^2)-d^2}$ & $e$ & $f$ \\
$a$ & $d$ & $e$ & $b$ & $c$ & $f$ \\
$a$ & $d$ & $e$ & $b$ & $c$ & $\sqrt{2(c^2+e^2)-f^2}$ \\
$a$ & $\sqrt{2(a^2+e^2)-d^2}$ & $e$ & $\sqrt{2(a^2+c^2)-b^2}$ & $c$ & $\sqrt{2(c^2+e^2)-f^2}$ \\
$a$ & $\sqrt{2(a^2+e^2)-d^2}$ & $e$ & $\sqrt{2(a^2+c^2)-b^2}$ & $c$ & $f$ \\
$c$ & $b$ & $a$ & $f$ & $e$ & $d$ \\
$c$ & $b$ & $a$ & $\sqrt{2(c^2+e^2)-f^2}$ & $e$ & $\sqrt{2(a^2+e^2)-d^2}$ \\
$c$ & $\sqrt{2(a^2+c^2)-b^2}$ & $a$ & $f$ & $e$ & $\sqrt{2(a^2+e^2)-d^2}$ \\
$c$ & $\sqrt{2(a^2+c^2)-b^2}$ & $a$ & $\sqrt{2(c^2+e^2)-f^2}$ & $e$ & $d$ \\
$c$ & $f$ & $e$ & $b$ & $a$ & $d$ \\
$c$ & $f$ & $e$ & $\sqrt{2(a^2+c^2)-b^2}$ & $a$ & $\sqrt{2(a^2+e^2)-d^2}$ \\
$c$ & $\sqrt{2(c^2+e^2)-f^2}$ & $e$ & $b$ & $a$ & $\sqrt{2(a^2+e^2)-d^2}$ \\
$c$ & $\sqrt{2(c^2+e^2)-f^2}$ & $e$ & $\sqrt{2(a^2+c^2)-b^2}$ & $a$ & $d$ \\
$e$ & $d$ & $a$ & $f$ & $c$ & $b$ \\
$e$ & $d$ & $a$ & $\sqrt{2(c^2+e^2)-f^2}$ & $c$ & $\sqrt{2(a^2+c^2)-b^2}$ \\
$e$ & $\sqrt{2(a^2+e^2)-d^2}$ & $a$ & $f$ & $c$ & $\sqrt{2(a^2+c^2)-b^2}$ \\
$e$ & $\sqrt{2(a^2+e^2)-d^2}$ & $a$ & $\sqrt{2(c^2+e^2)-f^2}$ & $c$ & $b$ \\
$e$ & $f$ & $c$ & $d$ & $a$ & $b$ \\
$e$ & $f$ & $c$ & $\sqrt{2(a^2+e^2)-d^2}$ & $a$ & $\sqrt{2(a^2+c^2)-b^2}$ \\
$e$ & $\sqrt{2(c^2+e^2)-f^2}$ & $c$ & $d$ & $a$ & $\sqrt{2(a^2+c^2)-b^2}$ \\
$e$ & $\sqrt{2(c^2+e^2)-f^2}$ & $c$ & $\sqrt{2(a^2+e^2)-d^2}$ & $a$ & $b$
\label{tbl:tetrapermute}
\end{longtable}
\end{center}

\vspace{-3em}
NOTE: While 24 tetrahedrons are usually in the group, sometimes less can occur,
due to symmetries of identical tetrahedron sides, which can create identical rows.
For {\it rectangular} {\it pipeds}, there are 6 tetrahedrons in the group.

\subsection{Density of rational solutions}

The next thing which must be considered is the density of the rational solutions in $\mathbb{Q}$ as compared
to the irrational solutions in $\mathbb{R}$ real space for a Diophantine geometrical object, when considering
the set of desired rational components. Sometimes this density is very sparse and clever algorithms must be utilized
to even find the rational solutions that match the desired set.

For example, while running computer studies of random parallelepipeds, it became necessary to create all possible parallelepipeds
using 2 integer sides and a 3rd integer diagonal, thus satisfying case 2, and case 5, for the parallelogram, automatically.

A record was kept of one such run using 3 integers, $a,b,c$, two $(a,b)$ for the 4 sides and $c$ for 1 rational
diagonal and ran the integers for $0 < a <= 100 \text{ and } 0 < b < a \text{ and } a-b < c < a+b$ to create
the parallelograms.
\begin{table}[ht]
\centering
\caption{Statistics for a short run of 746,344 integer sided parallelograms.}
\begin{tabular}{ r | c | c | l l}
count & percent & case & rational diagonal ? & rational area ? \\
\hline
737628 & 98.832\% & case 2 & no & no \\
1827 & 0.2448\% & case 5 & no & yes \\
6683 & 0.8954\% & case 3 & yes & no \\
206 & 0.0276\% & case 6 & yes & yes \\
\hline
63 & 0.00844\% & ibid & yes & right triangle \\
143 & 0.01916\% & ibid & yes & scalene triangle \\
\hline
746344 & 100.000\% & \textendash & \textendash &
\end{tabular}
\label{table:stats}
\end{table}

So it can be seen that rational solutions are sparse, only 206 solutions were found from 746,344 examined,
even for a computer algorithm which automatically satisfied four rational sides and 1 rational diagonal
$\in \mathbb{Z}$.

If we can parameterize the solutions and show that they completely cover the
rational solution space in $\mathbb{Z}$ or $\mathbb{Q}$, this is a vast improvement in
hunting for complete solutions for our set of properties to satisfy.

\subsection{Parameterizing Heron Triangles and Rational Parallelograms}\label{sub:params}

In order to efficiently discover rational solutions to \textit{Diophantine parallelepipeds}, it is very convenient to utilize
parametric solutions to automatically speed up the search process efficiently.

Such is the case here for both triangles and parallelograms.

The Heron Triangle is efficiently parameterized for the computer as the following solution using
4 integer parameters $m$, $n$, $p$, and $q$, all $\in \mathbb{Z}$:

\textbf{A parametric solution\cite{schubert} for Heron Triangles}
\begin{equation}
\begin{aligned}
\text{sides } \; a &= mn(p^2 + q^2) \\
b &= pq(m^2 + n^2) \\
c &= pq(n^2 - m^2) + mn(q^2 - p^2) \\
\triangle \text{ area} &= 4mnpq(mq+np)(nq-mp)
\end{aligned}
\label{eq:heron}
\end{equation}
This solution has been proven to fully cover the rational space $\mathbb{Q}$.

Another very helpful solution is that for parallelograms with 2 rational diagonals. This has been recently parameterized
by Walter Wyss\cite{wyss1}\cite{wyss2}. I changed his rational solution of $u,m,n \in \mathbb{Q}$ to that of
$k,m,n,p,q \in \mathbb{Z}$ since the integer case is easier to handle in computers than rational numbers.

In this parametric solution, the 5 integer parameters are $k$, $m$, $n$, $p$, and $q$. The scaling $k$ parameter creates
composite solutions(normally $k=1$). The 4 sides of the parallelogram are $a,b$ and the 2 diagonals are $c,d$.

\textbf{A parametric solution\cite{wyss2} for rational parallelograms}
\begin{equation}
\begin{aligned}
\text{    sides } \quad a & = k(nq - mp) \\
b & = k(mq + np) \\
\text{diagonals } \quad c & = k(p(m - n) + q(m + n)) \\
d & = k(p(n + m) + q(n - m))
\end{aligned}
\label{eq:para}
\end{equation}
This solution also fully covers rational space $\mathbb{Q}$

These 2 parameterization enormously speeded up the recovery of rational triangles and parallelograms used to assemble rational
component parallelepipeds, since the natural construction sequence is triangles $\rightarrow$ parallelograms $\rightarrow$ parallelepipeds.

\subsection{Density of Rational Solutions using Parameterization}

Even with parameterization, the density of rational solutions still has sparsity. In a run using
a parallelogram parameterization[\S\ref{sub:params} (\ref{eq:para})], satisfying case 3
automatically, 95,974,602 parallelograms with both diagonals rational were created. Only
25,088,615 were unique. Of those 25+ million, only 21,755 had rational area. This indicates that
in the search range of $0 < side < 10,001$ for the sides, only 0.0867\% satisfied the additional
rational constraint for case 6 using the parameterization.

Using the Heron parameterization[\S\ref{sub:params} (\ref{eq:heron})] which satisfies case 5 automatically,
of the 225,523 unique solutions found, with $0 < sides < 10,001$ only 5,302 had the 4th diagonal rational,
while 220,221 were irrational. This meant that only 2.351\% satisfied the constraint for case 6.
This is better than the previous 0.0867\% for the parallelogram 2 diagonals rational parameterization,
but 2\% is still low.

In another run of a search program using the Heron parametrization, 72,329,230 Heron triangles
were found but only 704,953 were unique. Of those, 686,264 had one diagonal irrational, only 18,689
satisfied case 6 with both diagonals rational and area rational. This means that only
2.651\% of the unique solutions were fully rational, but only 0.975\% of the initial solutions
satisfied case 6 constraints, for a net result of 0.0258\% of the raw solutions, even using
an efficient Heron triangle parameterization.

From this, it is seen that even when using parameterizations, the sparsity of rational solutions
satisfying constraints is low. It is necessary to use clever algorithms in computer searching
for Diophantine objects in embedded in 2d space. The requirements are more severe with objects
embedded in 3d space, and these constraints increase geometrically as we shall see.

\section{Computer Searches of Rational Parallelepipeds}

We offer some interesting discoveries of pipeds from the computer seaches.

To start, the extensive computer runs created 1,981,336,681 unique integer tetrahedrons which required over 441.8 gigabytes of storage space
to classify. The computer search also found 79,580 degenerate parallelepipeds, where the volume = 0, the algorithms creating the pipeds did not
specifically check for non-zero volume while assembling possible pipeds as this would have significantly increased the search time.

First of all, the following counts for the five classes of Diophantine parallelepipeds was obtained:\\
\begin{table}[h]
\begin{center}
\begin{tabular}{c|r|S[table-number-alignment=left] }
class & count \hspace{0.8 em} & \hphantom{00} \si{ \% \hspace{0.5em} abundance} \\
\hline
{\it acute triclinic} & 1315235647 & 66.3812293798 \\
{\it obtuse triclinic} & 659297660 & 33.2753976809 \\
{\it 1-ortho biclinic} & 6743400 & 0.340345992918 \\
{\it 2-ortho monoclinic} & 59930 & 0.00302472571041 \\
{\it rectangular} & 44 & 0.000002220723 \\
\hline
{\it Total count} & 1981336681 & 100.00 \\[0.5em]
\end{tabular}
\caption{Counts of Diophantine pipeds in each class}
\end{center}
\label{table:class-stats}
\end{table}

\vspace{-2em}
This shows that about 1 {\it obtuse} piped occurred for every 2 {\it acute} pipeds, while both of them together are $292.8\times$ as numerous
as {\it 1-ortho} pipeds. Both {\it acute} and {\it obtuse} pipeds together are $32,947.3\times$ more numerous than {\it 2-ortho} pipeds,
while {\it rectangular} pipeds are very scarce, accounting for only 0.0000022207\% of the tetrahedrons discovered. Futhermore
the ratio between the {\it 1-ortho} and {\it 2-ortho} pipeds is about 113 to 1.

There were 33,516 pipeds found which had rational volume, this is 0.0016915852\% or only about 1 in 59,116 pipeds.

\subsection{115 Unique Categories of Diophantine Pipeds}

After sorting the 1,981,336,681 integer tetrahedrons, it was discovered that 1,923 unique entries existed, grouped by the
27 rationality checks, and 115 unique categories resulted.

In the category table given below, we refer back to the 27 rationality checks previously determined in \S4.9 but we only supply the counts
in each rational check. For examine the skew triangle can have 16 possibilities, binary 4-place, decimal values 0...15, but we only give the total count of the rational components, if they occur, thus for the skew column, the count of rational occurances would be 0...4.

\noindent NOTE: For the volume check, -1 denotes 0 volume, no 3d object exists, while 0 denotes irrational volume, and 1 denotes
rational volume.

\begin{center}
\begin{longtable}{r|r|c|c|c|c|c|c|c|l}
\caption{115 Unique categories of Diophantine Parallelepipeds}\\
\hline
cat. & count & edges & skew & face & body & f area & b area & vol. & notes \\
\hline
\endfirsthead
\multicolumn{10}{c}
{\tablename\ \thetable\ -- \textit{115 Unique Categories - continued from previous page}} \\
\hline
cat. & count & edges & skew & face & body & f area & b area & vol. & notes \\
\hline
\endhead
\hline \multicolumn{10}{r}{\textit{Continued on next page}} \\
\endfoot
\hline
\endlastfoot
1 & 268 & 3 & 0 & 6 & 0 & 0 & 0 & -1 & \\
2 & 1914630558 & 3 & 0 & 6 & 0 & 0 & 0 & 0 & \\
3 & 14628 & 3 & 0 & 6 & 0 & 0 & 0 & 1 & \\
4 & 2585672 & 3 & 0 & 6 & 0 & 0 & 1 & 0 & \\
5 & 48 & 3 & 0 & 6 & 0 & 0 & 1 & 1 & \\
6 & 13252 & 3 & 0 & 6 & 0 & 0 & 2 & 0 & \\
7 & 96 & 3 & 0 & 6 & 0 & 0 & 2 & 1 & \\
8 & 48 & 3 & 0 & 6 & 0 & 0 & 3 & 0 & \\
9 & 32448 & 3 & 0 & 6 & 0 & 1 & 0 & -1 & \\
10 & 19497256 & 3 & 0 & 6 & 0 & 1 & 0 & 0 & \\
11 & 912 & 3 & 0 & 6 & 0 & 1 & 0 & 1 & \\
12 & 6831972 & 3 & 0 & 6 & 0 & 1 & 1 & 0 & \\
13 & 72 & 3 & 0 & 6 & 0 & 1 & 1 & 1 & \\
14 & 1412 & 3 & 0 & 6 & 0 & 1 & 2 & 0 & \\
15 & 2304 & 3 & 0 & 6 & 0 & 1 & 2 & 1 & \\
16 & 12232 & 3 & 0 & 6 & 0 & 1 & 3 & 0 & \\
17 & 69136 & 3 & 0 & 6 & 0 & 2 & 0 & 0 & \\
18 & 2368 & 3 & 0 & 6 & 0 & 2 & 1 & 0 & \\
19 & 20898 & 3 & 0 & 6 & 0 & 2 & 2 & 0 & \\
20 & 110 & 3 & 0 & 6 & 0 & 2 & 4 & 0 & \\
21 & 4440 & 3 & 0 & 6 & 0 & 3 & 0 & 0 & \\
22 & 19628 & 3 & 0 & 6 & 0 & 3 & 1 & 0 & \\
23 & 7962 & 3 & 0 & 6 & 0 & 3 & 2 & 1 & \\
24 & 3776 & 3 & 0 & 6 & 0 & 3 & 3 & 0 & \\
25 & 44 & 3 & 0 & 6 & 0 & 3 & 6 & 1 & {\it rect} \\
26 & 24 & 3 & 0 & 6 & 1 & 0 & 0 & -1 & \\
27 & 32314073 & 3 & 0 & 6 & 1 & 0 & 0 & 0 & \\
28 & 5036 & 3 & 0 & 6 & 1 & 0 & 0 & 1 & \\
29 & 145920 & 3 & 0 & 6 & 1 & 0 & 1 & 0 & \\
30 & 192 & 3 & 0 & 6 & 1 & 0 & 1 & 1 & \\
31 & 2648 & 3 & 0 & 6 & 1 & 0 & 2 & 0 & \\
32 & 5232 & 3 & 0 & 6 & 1 & 1 & 0 & -1 & \\
33 & 477828 & 3 & 0 & 6 & 1 & 1 & 0 & 0 & \\
34 & 384 & 3 & 0 & 6 & 1 & 1 & 0 & 1 & \\
35 & 65516 & 3 & 0 & 6 & 1 & 1 & 1 & 0 & \\
36 & 24 & 3 & 0 & 6 & 1 & 1 & 2 & 0 & \\
37 & 208 & 3 & 0 & 6 & 1 & 1 & 2 & 1 & \\
38 & 196 & 3 & 0 & 6 & 1 & 1 & 3 & 0 & \\
39 & 1892 & 3 & 0 & 6 & 1 & 2 & 0 & 0 & \\
40 & 216 & 3 & 0 & 6 & 1 & 2 & 1 & 0 & \\
41 & 264 & 3 & 0 & 6 & 1 & 3 & 1 & 0 & \\
42 & 834952 & 3 & 0 & 6 & 2 & 0 & 0 & 0 & \\
43 & 264 & 3 & 0 & 6 & 2 & 0 & 0 & 1 & \\
44 & 24 & 3 & 0 & 6 & 2 & 0 & 1 & -1 & \\
45 & 12240 & 3 & 0 & 6 & 2 & 0 & 1 & 0 & \\
46 & 504 & 3 & 0 & 6 & 2 & 0 & 2 & 0 & \\
47 & 3056 & 3 & 0 & 6 & 2 & 1 & 0 & -1 & \\
48 & 102768 & 3 & 0 & 6 & 2 & 1 & 0 & 0 & \\
49 & 53252 & 3 & 0 & 6 & 2 & 1 & 1 & 0 & \\
50 & 1320 & 3 & 0 & 6 & 2 & 1 & 3 & 0 & \\
51 & 24 & 3 & 0 & 6 & 2 & 2 & 0 & 0 & \\
52 & 718 & 3 & 0 & 6 & 2 & 2 & 2 & 0 & \\
53 & 72 & 3 & 0 & 6 & 2 & 3 & 1 & 0 & \\
54 & 182 & 3 & 0 & 6 & 2 & 3 & 2 & 1 & \\
55 & 19416 & 3 & 0 & 6 & 3 & 0 & 0 & 0 & \\
56 & 192 & 3 & 0 & 6 & 3 & 0 & 1 & 0 & \\
57 & 64 & 3 & 0 & 6 & 3 & 1 & 0 & -1 & \\
58 & 5160 & 3 & 0 & 6 & 3 & 1 & 0 & 0 & \\
59 & 1912 & 3 & 0 & 6 & 3 & 1 & 1 & 0 & \\
60 & 124 & 3 & 0 & 6 & 4 & 0 & 0 & 0 & {\it perfect} \\
61 & 600 & 3 & 0 & 6 & 4 & 0 & 3 & -1 & \\
62 & 4 & 3 & 0 & 6 & 4 & 1 & 0 & 0 & {\it perfect} \\
63 & 1324 & 3 & 1 & 6 & 0 & 0 & 0 & -1 & \\
64 & 3126152 & 3 & 1 & 6 & 0 & 0 & 0 & 0 & \\
65 & 96 & 3 & 1 & 6 & 0 & 0 & 0 & 1 & \\
66 & 29040 & 3 & 1 & 6 & 0 & 0 & 1 & 0 & \\
67 & 648 & 3 & 1 & 6 & 0 & 0 & 2 & 0 & \\
68 & 74084 & 3 & 1 & 6 & 0 & 1 & 0 & 0 & \\
69 & 336 & 3 & 1 & 6 & 0 & 1 & 0 & 1 & \\
70 & 83088 & 3 & 1 & 6 & 0 & 1 & 1 & 0 & \\
71 & 144 & 3 & 1 & 6 & 0 & 1 & 1 & 1 & \\
72 & 24 & 3 & 1 & 6 & 0 & 1 & 2 & 0 & \\
73 & 248 & 3 & 1 & 6 & 0 & 1 & 2 & 1 & \\
74 & 300 & 3 & 1 & 6 & 0 & 1 & 3 & 0 & \\
75 & 312 & 3 & 1 & 6 & 0 & 2 & 0 & 0 & \\
76 & 108 & 3 & 1 & 6 & 0 & 3 & 1 & 0 & \\
77 & 120 & 3 & 1 & 6 & 1 & 0 & 0 & -1 & \\
78 & 128176 & 3 & 1 & 6 & 1 & 0 & 0 & 0 & \\
79 & 2988 & 3 & 1 & 6 & 1 & 0 & 1 & 0 & \\
80 & 5160 & 3 & 1 & 6 & 1 & 1 & 0 & 0 & \\
81 & 24 & 3 & 1 & 6 & 1 & 1 & 0 & 1 & \\
82 & 2880 & 3 & 1 & 6 & 1 & 1 & 1 & 0 & \\
83 & 24 & 3 & 1 & 6 & 1 & 1 & 2 & 1 & \\
84 & 72 & 3 & 1 & 6 & 1 & 2 & 0 & 0 & \\
85 & 5136 & 3 & 1 & 6 & 2 & 0 & 0 & 0 & \\
86 & 48 & 3 & 1 & 6 & 2 & 0 & 0 & 1 & \\
87 & 1304 & 3 & 1 & 6 & 2 & 0 & 1 & -1 & \\
88 & 360 & 3 & 1 & 6 & 2 & 1 & 0 & 0 & \\
89 & 2460 & 3 & 1 & 6 & 2 & 1 & 1 & 0 & \\
90 & 432 & 3 & 1 & 6 & 2 & 1 & 2 & 0 & \\
91 & 192 & 3 & 1 & 6 & 2 & 1 & 3 & 0 & \\
92 & 72 & 3 & 1 & 6 & 3 & 0 & 0 & 0 & \\
93 & 12 & 3 & 1 & 6 & 3 & 1 & 1 & 0 & \\
94 & 16352 & 3 & 2 & 6 & 0 & 0 & 0 & 0 & \\
95 & 24 & 3 & 2 & 6 & 0 & 0 & 0 & 1 & \\
96 & 916 & 3 & 2 & 6 & 0 & 0 & 1 & 0 & \\
97 & 1476 & 3 & 2 & 6 & 0 & 1 & 0 & 0 & \\
98 & 3388 & 3 & 2 & 6 & 0 & 1 & 1 & 0 & \\
99 & 36 & 3 & 2 & 6 & 0 & 1 & 2 & 1 & \\
100 & 12 & 3 & 2 & 6 & 0 & 1 & 3 & 0 & \\
101 & 96 & 3 & 2 & 6 & 0 & 2 & 0 & 0 & \\
102 & 12 & 3 & 2 & 6 & 0 & 2 & 2 & 0 & \\
103 & 36 & 3 & 2 & 6 & 0 & 3 & 1 & 0 & \\
104 & 12 & 3 & 2 & 6 & 0 & 3 & 2 & 1 & \\
105 & 808 & 3 & 2 & 6 & 1 & 0 & 0 & 0 & \\
106 & 3024 & 3 & 2 & 6 & 1 & 1 & 0 & 0 & \\
107 & 1404 & 3 & 2 & 6 & 1 & 1 & 1 & 0 & \\
108 & 192 & 3 & 2 & 6 & 1 & 1 & 2 & 1 & \\
109 & 48 & 3 & 2 & 6 & 2 & 1 & 0 & 0 & \\
110 & 144 & 3 & 2 & 6 & 2 & 1 & 1 & 0 & \\
111 & 180 & 3 & 3 & 6 & 0 & 1 & 1 & 0 & \\
112 & 32744 & 3 & 4 & 6 & 0 & 3 & 6 & -1 & \\
113 & 48 & 3 & 4 & 6 & 1 & 3 & 6 & -1 & \\
114 & 2300 & 3 & 4 & 6 & 2 & 3 & 6 & -1 & \\
115 & 24 & 3 & 4 & 6 & 4 & 3 & 6 & -1 & planar\\
\hline \multicolumn{10}{l}{\textit{1,981,336,681 Diophantine pipeds in 1923 entries in 115 categories}}
\label{tbl:115}
\end{longtable}
\end{center}

We will select interesting examples from these 115 categories of pipeds in the next sections.

\subsection{Perfect Diophantine parallelepipeds}

\textit{Perfect Diophantine parallelepipeds} are those where both the 6 face diagonals and the 4 body diagonals
are rational.

The first perfect parallelepiped discovered was the Diophantine \textit{acute piped} created by these 3 vectors:
\begin{center}
\begin{tabular}{r|cr}
vector & direction vector & length \\
\hline
$v_1$ & $[106, 0, 0]$ & 106 \\
$v_2$ & $[\frac{4913}{53}, \frac{30\sqrt{202398}}{53}, 0]$ & 271 \\
$v_3$ & $[\frac{2911}{53}, \frac{3468\sqrt{202398}}{255407}, \frac{66\sqrt{40277202}}{4819}]$ & 103 \\
\end{tabular}
\end{center}

\newpage

\begin{figure}[!ht]
\centering
\includegraphics[scale=1.0]{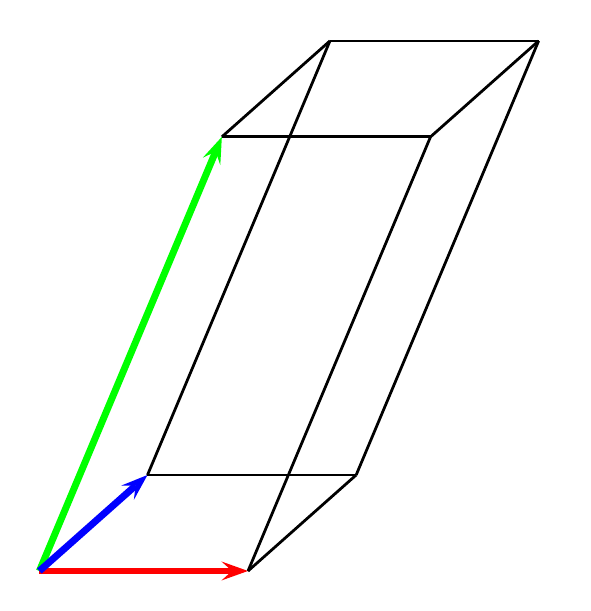}
\captionsetup{justification=centering}
\caption[caption]{The Sawyer-Reiter Acute Diophantine Piped \\ \vspace{-10pt} \hspace{\textwidth} {\small $\alpha = \arccos\left(\frac{4913}{14363}\right) \; \beta = \arccos\left(\frac{2911}{5459}\right) \; \gamma = \arccos\left(\frac{6647}{27913}\right) \text{ and } a \ne b \ne c $ }}
\label{fig:sawyer-reiter}
\end{figure}

This \textit{acute} piped discovered by Jorge F. Sawyer and Clifford A. Reiter\cite{sawyer} in 2008 or
2009 was the first known that had all 4 body diagonals rational as well as all 12 face diagonals rational.

\noindent
The computer searches found 6 \textit{perfect Diophantine parallelepipeds} including the Sawyer-Reiter perfect piped shown
in Fig. \ref{fig:sawyer-reiter}.
\begin{table}[htpb]
\centering
\caption{Six perfect Diophantine parallelepipeds}
{\footnotesize
\begin{tabular}{| c c c | c c c c | c c c c c c | c c c c | c c c | c c c c c c | c |}
\multicolumn{3}{c}{edges} &
\multicolumn{4}{c}{skew triangles} &
\multicolumn{6}{c}{face diagonals} &
\multicolumn{4}{c}{body diagonals} &
\multicolumn{3}{c}{face area} &
\multicolumn{6}{c}{body area} &
\multicolumn{1}{c}{vol} \\
\hline
 1 & 1 & 1 & 0 & 0 & 0 & 0 & 1 & 1 & 1 & 1 & 1 & 1 & 1 & 1 & 1 & 1 & 0 & 0 & 0 & 0 & 0 & 0 & 0 & 0 & 0 & 0 \Tstrut \\
\hline
\end{tabular}
}
{\tiny
\begin{tabular}{ r | c c c | c c c c c c }
 & \multicolumn{3}{c}{Basis Vectors} & \multicolumn{6}{c}{Tetrahedron Sides} \\
Class & $v_1$ & $v_2$ & $v_3$ & a & b & c & d & e & f \\[0.15em]
\hline\\[-0.5em]
acute & $\left(103,0,0\right)$ & $\left(\frac{5822}{103},\frac{84\sqrt{12090}}{103},0\right)$ & $\left(\frac{6647}{103},\frac{26010\sqrt{12090}}{41509},\frac{66\sqrt{2405910}}{403}\right)$ & 103 & 101 & 106 & 266 & 271 & 255 \Tstrut \Bstrut \\[0.5em]
acute & $\left(335,0,0\right)$ & $\left(\frac{24157}{335},\frac{96\sqrt{1558986}}{335},0\right)$ & $\left(\frac{15510}{67},\frac{2089315\sqrt{1558986}}{5356516},\frac{65\sqrt{160551933366}}{79948}\right)$ & 335 & 444 & 365 & 595 & 630 & 385 \Tstrut \Bstrut \\[0.5em]
1-ortho & $\left(340,0,0\right)$ & $\left(0,357,0\right)$ & $\left(\frac{7400}{17},\frac{6384}{17},\frac{720\sqrt{321}}{17}\right)$ & 340 & 493 & 357 & 852 & 952 & 875 \Tstrut \Bstrut \\[0.5em]
acute & $\left(342,0,0\right)$ & $\left(\frac{21385}{57},\frac{70\sqrt{141410}}{57},0\right)$ & $\left(\frac{18847}{57},\frac{614608\sqrt{141410}}{4030185},\frac{18\sqrt{6688793542510}}{70705}\right)$ & 342 & 463 & 595 & 661 & 739 & 774 \Tstrut \Bstrut \\[0.5em]
acute & $\left(375,0,0\right)$ & $\left(468,72\sqrt{14},0\right)$ & $\left(\frac{11951}{25},-\frac{6088\sqrt{14}}{75},\frac{8\sqrt{343966}}{15}\right)$ & 375 & 285 & 540 & 448 & 647 & 653 \Tstrut \Bstrut \\[0.5em]
acute & $\left(422,0,0\right)$ & $\left(\frac{81891}{211},\frac{42\sqrt{4659370}}{211},0\right)$ & $\left(\frac{175195}{211},-\frac{2958000\sqrt{4659370}}{98312707},\frac{30\sqrt{39076439247446}}{465937}\right)$ & 422 & 431 & 579 & 577 & 925 & 776 \Tstrut \Bstrut \\[0.5em]
\hline
\end{tabular}
}
\label{table:sixperf}
\end{table}
Further computer searches have revealed that the {\it 2-ortho} {\it monoclinic} {\it piped} can have all diagonals rational also. Rathbun found parametric formuli for such pipeds\cite{rathbun1},\cite{rathbun2}.

\newpage

\subsection{Interesting examples of Diophantine pipeds found by searches}

In the 14 examples given below, in general, the first occurrences of the designated piped is given.

\subsubsection{Rational volume pipeds}
Searching found that Diophantine parallelepipeds can have a rational volume, although none of the body diagonals have a rational
length, nor does any of the 6 body or 3 face parallelograms have a rational area.
\begin{table}[h]
\centering
\caption{Rational volume Diophantine parallelepipeds}
{\scriptsize
\begin{tabular}{| c c c | c c c c | c c c c c c | c c c c | c c c | c c c c c c | c |}
\multicolumn{3}{c}{edges} &
\multicolumn{4}{c}{skew triangles} &
\multicolumn{6}{c}{face diagonals} &
\multicolumn{4}{c}{body diagonals} &
\multicolumn{3}{c}{face area} &
\multicolumn{6}{c}{body area} &
\multicolumn{1}{c}{vol} \\
\hline
 1 & 1 & 1 & 0 & 0 & 0 & 0 & 1 & 1 & 1 & 1 & 1 & 1 & 0 & 0 & 0 & 0 & 0 & 0 & 0 & 0 & 0 & 0 & 0 & 0 & 0 & 1 \\
\hline
\end{tabular}
\begin{tabular}{ r | c c c | c c c c c c | c }
& \multicolumn{3}{c}{Basis Vectors} & \multicolumn{6}{|c}{Tetrahedron Sides} & \multicolumn{1}{|c}{} \\
Class & $v_1$ & $v_2$ & $v_3$ & a & b & c & d & e & f & Volume \\[0.15em]
\hline\\[-0.5em]
acute & $\left(17,0,0\right)$ & $\left(\frac{473}{17},\frac{24\sqrt{455}}{17},0\right)$ & $\left(\frac{876}{17},\frac{1848\sqrt{455}}{1105},\frac{108\sqrt{455}}{65}\right)$ & 17 & 32 & 41 & 61 & 72 & 43 & 18144 \Tstrut \Bstrut \\[0.5em]
obtuse & $\left(19,0,0\right)$ & $\left(\frac{122}{19},\frac{84\sqrt{170}}{19},0\right)$ & $\left(\frac{203}{19},-\frac{216\sqrt{170}}{323},\frac{24\sqrt{170}}{17}\right)$ & 19 & 59 & 58 & 22 & 23 & 69 & 20160 \Tstrut \Bstrut \\[0.5em]
acute & $\left(25,0,0\right)$ & $\left(\frac{342}{25},\frac{12\sqrt{3094}}{25},0\right)$ & $\left(22,\frac{132\sqrt{3094}}{221},\frac{120\sqrt{3094}}{221}\right)$ & 25 & 29 & 30 & 45 & 50 & 32 & 20160 \Tstrut \Bstrut \\[0.5em]
acute & $\left(26,0,0\right)$ & $\left(\frac{207}{13},\frac{30\sqrt{238}}{13},0\right)$ & $\left(\frac{216}{13},\frac{444\sqrt{238}}{221},\frac{36\sqrt{238}}{17}\right)$ & 26 & 37 & 39 & 46 & 48 & 33 & 30240 \Tstrut \Bstrut \\[0.5em]
acute & $\left(28,0,0\right)$ & $\left(14,21\sqrt{5},0\right)$ & $\left(18,\frac{111\sqrt{5}}{5},\frac{48\sqrt{5}}{5}\right)$ & 28 & 49 & 49 & 55 & 57 & 22 & 28224 \Tstrut \Bstrut \\[0.5em]
\hline
\end{tabular}
}
\label{table:5ratvol}
\end{table}
Interesting enough, the $2^{\text nd}$ and the $3^{\text rd}$ pipeds in Table \ref{table:5ratvol}. have the same volume, 20160.

\subsubsection{Rational pipeds with one body parallelogram area rational and rational volume}
Diophantine pipeds were found that had one body parallelogram with rational area, while the other 5 are irrational. None of
the face parallelograms had rational area. The volume was rational also.
\begin{table}[h]
\centering
\caption{Rational pipeds with one rational body parallelogram area and rational volume}
{\scriptsize
\begin{tabular}{ r | c c c | c c c c c c }
 & \multicolumn{3}{c}{Basis Vectors} & \multicolumn{6}{c}{Tetrahedron Sides} \\
Class & $v_1$ & $v_2$ & $v_3$ & a & b & c & d & e & f \\[0.15em]
\hline\\[-0.5em]
obtuse & $\left(99,0,0\right)$ & $\left(\frac{301}{3},\frac{104\sqrt{170}}{3},0\right)$ & $\left(\frac{77}{3},-\frac{17072\sqrt{170}}{3315},\frac{16632\sqrt{170}}{1105}\right)$ & 99 & 452 & 463 & 220 & 209 & 560 \Tstrut \Bstrut \\[0.5em]
obtuse & $\left(99,0,0\right)$ & $\left(-\frac{301}{3},\frac{104\sqrt{170}}{3},0\right)$ & $\left(-\frac{77}{3},-\frac{17072\sqrt{170}}{3315},\frac{16632\sqrt{170}}{1105}\right)$ & 99 & 494 & 463 & 242 & 209 & 560 \Tstrut \Bstrut \\[0.5em]
acute & $\left(274,0,0\right)$ & $\left(\frac{42364}{137},\frac{120\sqrt{114130}}{137},0\right)$ & $\left(\frac{49679}{137},\frac{1829856\sqrt{114130}}{1563581},\frac{10296\sqrt{114130}}{11413}\right)$ & 274 & 298 & 428 & 507 & 617 & 325 \Tstrut \Bstrut \\[0.5em]
obtuse & $\left(209,0,0\right)$ & $\left(\frac{2525}{19},\frac{336\sqrt{629}}{19},0\right)$ & $\left(\frac{231}{19},-\frac{12408\sqrt{629}}{11951},\frac{2376\sqrt{629}}{629}\right)$ & 209 & 450 & 463 & 220 & 99 & 494 \Tstrut \Bstrut \\[0.5em]
obtuse & $\left(209,0,0\right)$ & $\left(\frac{2525}{19},\frac{336\sqrt{629}}{19},0\right)$ & $\left(-\frac{231}{19},\frac{12408\sqrt{629}}{11951},\frac{2376\sqrt{629}}{629}\right)$ & 209 & 450 & 463 & 242 & 99 & 452 \Tstrut \Bstrut \\[0.5em]
\hline
\end{tabular}
}
\label{table:bodyarea}
\end{table}
\begin{table}[htpb]
\centering
\caption{Rational component checks for the rational body area pipeds}
{\footnotesize
\begin{tabular}{| c c c | c c c c | c c c c c c | c c c c | c c c | c c c c c c | c |}
\multicolumn{3}{c}{edges} &
\multicolumn{4}{c}{skew triangles} &
\multicolumn{6}{c}{face diagonals} &
\multicolumn{4}{c}{body diagonals} &
\multicolumn{3}{c}{face area} &
\multicolumn{6}{c}{body area} &
\multicolumn{1}{c}{vol} \\
\hline
 1 & 1 & 1 & 0 & 0 & 0 & 0 & 1 & 1 & 1 & 1 & 1 & 1 & 0 & 0 & 0 & 0 & 0 & 0 & 0 & 1 & 0 & 0 & 0 & 0 & 0 & 1 \\
 1 & 1 & 1 & 0 & 0 & 0 & 0 & 1 & 1 & 1 & 1 & 1 & 1 & 0 & 0 & 0 & 0 & 0 & 0 & 0 & 1 & 0 & 0 & 0 & 0 & 0 & 1 \\
 1 & 1 & 1 & 0 & 0 & 0 & 0 & 1 & 1 & 1 & 1 & 1 & 1 & 0 & 0 & 0 & 0 & 0 & 0 & 0 & 0 & 0 & 0 & 1 & 0 & 0 & 1 \\
 1 & 1 & 1 & 0 & 0 & 0 & 0 & 1 & 1 & 1 & 1 & 1 & 1 & 0 & 0 & 0 & 0 & 0 & 0 & 0 & 0 & 0 & 0 & 0 & 0 & 1 & 1 \\
 1 & 1 & 1 & 0 & 0 & 0 & 0 & 1 & 1 & 1 & 1 & 1 & 1 & 0 & 0 & 0 & 0 & 0 & 0 & 0 & 0 & 0 & 0 & 0 & 0 & 1 & 1 \\
\hline
\end{tabular}
}
\label{table:bodyareasig}
\end{table}

The volume of these 5 pipeds are 8781696, 8781696, 24710400, 8781696, 8781696. As shown by these identical volumes
the $1^{\text{st}}$, $2^{\text{nd}}$, $4^{\text{th}}$, and $5^{\text{th}}$ tetrahedrons are in the same family, while the
 $3^{\text{rd}}$ is in another family.

\newpage

\subsubsection{Rational pipeds with all 9 parallelograms of rational area and rational volume}

It was discovered that the {\it rectangular} class of Diophantine pipeds , has all 3 face parallelograms and all 6 body parallelograms
with rational area, as well as the volume rational. Unfortunately none of the 4 body diagonals were rational. This class of piped
 is related to the question of whether or not {\it perfect} {\it cuboids} exist, still unanswered.

\begin{table}[h]
\centering
\caption{Parellelepipeds with all face and body parallelograms of rational area}
{\footnotesize
\begin{tabular}{| c c c | c c c c | c c c c c c | c c c c | c c c | c c c c c c | c |}
\multicolumn{3}{c}{edges} &
\multicolumn{4}{c}{skew triangles} &
\multicolumn{6}{c}{face diagonals} &
\multicolumn{4}{c}{body diagonals} &
\multicolumn{3}{c}{face area} &
\multicolumn{6}{c}{body area} &
\multicolumn{1}{c}{vol} \\
\hline
 1 & 1 & 1 & 0 & 0 & 0 & 0 & 1 & 1 & 1 & 1 & 1 & 1 & 0 & 0 & 0 & 0 & 1 & 1 & 1 & 1 & 1 & 1 & 1 & 1 & 1 & 1 \\
\hline
\end{tabular}
}
{\footnotesize
\begin{tabular}{ r | c c c | c c c c c c | l}
 & \multicolumn{3}{c|}{Basis Vectors} & \multicolumn{6}{c|}{Tetrahedron Sides} & \\
Class & $v_1$ & $v_2$ & $v_3$ & a & b & c & d & e & f & Volume \\[0.15em]
\hline\\[-0.5em]
rectangular & $\left(44,0,0\right)$ & $\left(0,117,0\right)$ & $\left(0,0,240\right)$ & 44 & 125 & 117 & 244 & 240 & 267 & 1235520 \Tstrut \Bstrut \\[0.5em]
rectangular & $\left(85,0,0\right)$ & $\left(0,132,0\right)$ & $\left(0,0,720\right)$ & 85 & 157 & 132 & 725 & 720 & 732 & 8078400 \Tstrut \Bstrut \\[0.5em]
rectangular & $\left(140,0,0\right)$ & $\left(0,480,0\right)$ & $\left(0,0,693\right)$ & 140 & 500 & 480 & 707 & 693 & 843 & 46569600 \Tstrut \Bstrut \\[0.5em]
rectangular & $\left(160,0,0\right)$ & $\left(0,231,0\right)$ & $\left(0,0,792\right)$ & 160 & 281 & 231 & 808 & 792 & 825 & 29272320 \Tstrut \Bstrut \\[0.5em]
rectangular & $\left(187,0,0\right)$ & $\left(0,1020,0\right)$ & $\left(0,0,1584\right)$ & 187 & 1037 & 1020 & 1595 & 1584 & 1884 & 302132160 \Tstrut \Bstrut \\[0.5em]
rectangular & $\left(240,0,0\right)$ & $\left(0,252,0\right)$ & $\left(0,0,275\right)$ & 240 & 348 & 252 & 365 & 275 & 373 & 16632000 \Tstrut \Bstrut \\[0.5em]
\hline
\end{tabular}
}
\label{table:3face6body}
\end{table}

\subsubsection{Rational piped with 1 rational body diagonal, and some rational area parallelograms}

Computer searches led to the discovery of one piped which had 1 body diagonal rational, and 2 body parallelograms
with rational area, as well as one face parallelogram with rational area.

\begin{table}[h]
\centering
\caption{Diophantine piped with a rational body diagonal and some parallelograms of rational area}
{\footnotesize
\begin{tabular}{| c c c | c c c c | c c c c c c | c c c c | c c c | c c c c c c | c |}
\multicolumn{3}{c}{edges} &
\multicolumn{4}{c}{skew triangles} &
\multicolumn{6}{c}{face diagonals} &
\multicolumn{4}{c}{body diagonals} &
\multicolumn{3}{c}{face area} &
\multicolumn{6}{c}{body area} &
\multicolumn{1}{c}{vol} \\
\hline
 1 & 1 & 1 & 0 & 0 & 0 & 0 & 1 & 1 & 1 & 1 & 1 & 1 & 0 & 0 & 0 & 1 & 1 & 0 & 0 & 1 & 0 & 0 & 0 & 0 & 1 & 0 \\
\hline
\end{tabular}
}
{\footnotesize
\begin{tabular}{ r | c c c | c c c c c c }
 & \multicolumn{3}{c|}{Basis Vectors} & \multicolumn{6}{c}{Tetrahedron Sides} \\
Class & $v_1$ & $v_2$ & $v_3$ & a & b & c & d & e & f \\[0.15em]
\hline\\[-0.5em]
acute & $\left(385,0,0\right)$ & $\left(\frac{1211}{5},\frac{24\sqrt{4929}}{5},0\right)$ & $\left(\frac{539}{5},\frac{14168\sqrt{4929}}{8215},\frac{2464\sqrt{54219}}{1643}\right)$ & 385 & 366 & 415 & 462 & 385 & 432 \Tstrut \Bstrut \\[0.5em]
\hline
\end{tabular}
}
\label{table:1bd1fa2bd}
\end{table}

\subsubsection{Rational piped with 2 rational body diagonals and 2 rational area face parallelograms}
One Diophantine parallelepiped was found that had 2 rational body diagonals (of 4 possible) and 2 rational area
face parallelograms (of 3 possible).
\begin{table}[h]
\centering
\caption{Diophantine piped with 2 rational body diagonals and 2 rational area face parallelograms}
{\footnotesize
\begin{tabular}{| c c c | c c c c | c c c c c c | c c c c | c c c | c c c c c c | c |}
\multicolumn{3}{c}{edges} &
\multicolumn{4}{c}{skew triangles} &
\multicolumn{6}{c}{face diagonals} &
\multicolumn{4}{c}{body diagonals} &
\multicolumn{3}{c}{face area} &
\multicolumn{6}{c}{body area} &
\multicolumn{1}{c}{vol} \\
\hline
 1 & 1 & 1 & 0 & 0 & 0 & 0 & 1 & 1 & 1 & 1 & 1 & 1 & 0 & 1 & 0 & 1 & 1 & 0 & 1 & 0 & 0 & 0 & 0 & 0 & 0 & 0 \\
\hline
\end{tabular}
}
{\footnotesize
\begin{tabular}{ r | c c c | c c c c c c }
 & \multicolumn{3}{c|}{Basis Vectors} & \multicolumn{6}{c}{Tetrahedron Sides} \\
Class & $v_1$ & $v_2$ & $v_3$ & a & b & c & d & e & f \\[0.15em]
\hline\\[-0.5em]
1-ortho & $\left(175,0,0\right)$ & $\left(\frac{4017}{175},\frac{792\sqrt{3149}}{175},0\right)$ & $\left(0,\frac{122500\sqrt{3149}}{34639},\frac{280\sqrt{2096947441}}{34639}\right)$ & 175 & 296 & 255 & 455 & 420 & 375 \Tstrut \Bstrut \\[0.5em]
\hline
\end{tabular}
}
\label{table:2bd2fa}
\end{table}

\newpage

\subsubsection{Rational piped with 4 rational body diagonals and 1 rational area face parallelogram}
One Diophantine piped was found that had all 4 body diagonals rational, and had one face {\paragram } parallelogram
with rational area.
\begin{table}[h]
\centering
\caption{Diophantine piped with 4 rational body diagonals and 1 rational area face parallelogram }
{\footnotesize
\begin{tabular}{| c c c | c c c c | c c c c c c | c c c c | c c c | c c c c c c | c |}
\multicolumn{3}{c}{edges} &
\multicolumn{4}{c}{skew triangles} &
\multicolumn{6}{c}{face diagonals} &
\multicolumn{4}{c}{body diagonals} &
\multicolumn{3}{c}{face area} &
\multicolumn{6}{c}{body area} &
\multicolumn{1}{c}{vol} \\
\hline
 1 & 1 & 1 & 0 & 0 & 0 & 0 & 1 & 1 & 1 & 1 & 1 & 1 & 1 & 1 & 1 & 1 & 0 & 1 & 0 & 0 & 0 & 0 & 0 & 0 & 0 & 0 \\
\hline
\end{tabular}
}
{\footnotesize
\begin{tabular}{ r | c c c | c c c c c c }
 & \multicolumn{3}{c|}{Basis Vectors} & \multicolumn{6}{c}{Tetrahedron Sides} \\
Class & $v_1$ & $v_2$ & $v_3$ & a & b & c & d & e & f \\[0.15em]
\hline\\[-0.5em]
1-ortho & $\left(340,0,0\right)$ & $\left(0,357,0\right)$ & $\left(\frac{7400}{17},\frac{6384}{17},\frac{720\sqrt{321}}{17}\right)$ & 340 & 493 & 357 & 852 & 952 & 875 \Tstrut \Bstrut \\[0.5em]
\hline
\end{tabular}
}
\label{table:4bd1fa}
\end{table}

\subsubsection{Rational piped with 1 rational area face and 2 rational area body parallelograms}
One piped was found that had a rational face area parallelogram and 2 rational area body parallelograms.
\begin{table}[h]
\centering
\caption{Diophantine piped with 1 rational area face and 2 rational area body parallelograms}
{\footnotesize
\begin{tabular}{| c c c | c c c c | c c c c c c | c c c c | c c c | c c c c c c | c |}
\multicolumn{3}{c}{edges} &
\multicolumn{4}{c}{skew triangles} &
\multicolumn{6}{c}{face diagonals} &
\multicolumn{4}{c}{body diagonals} &
\multicolumn{3}{c}{face area} &
\multicolumn{6}{c}{body area} &
\multicolumn{1}{c}{vol} \\
\hline
 1 & 1 & 1 & 0 & 0 & 0 & 1 & 1 & 1 & 1 & 1 & 1 & 1 & 0 & 0 & 0 & 0 & 1 & 0 & 0 & 1 & 0 & 0 & 0 & 1 & 0 & 0 \\
\hline
\end{tabular}
}
{\footnotesize
\begin{tabular}{ r | c c c | c c c c c c }
 & \multicolumn{3}{c|}{Basis Vectors} & \multicolumn{6}{c}{Tetrahedron Sides} \\
Class & $v_1$ & $v_2$ & $v_3$ & a & b & c & d & e & f \\[0.15em]
\hline\\[-0.5em]
1-ortho & $\left(204,0,0\right)$ & $\left(\frac{2400}{17},\frac{720\sqrt{21}}{17},0\right)$ & $\left(0,\frac{680\sqrt{21}}{21},\frac{85\sqrt{20265}}{21}\right)$ & 204 & 204 & 240 & 629 & 595 & 595 \Tstrut \Bstrut \\[0.5em]
\hline
\end{tabular}
}
\label{table:1fa2ba}
\end{table}

\subsubsection{Rational piped with 1 rational body diagonal, 1 rational face parallelogram and 2 rational body parallelograms and rational
volume}
One Diophantine piped was found which has 1 rational body diagonal, 1 face parallelogram with rational area, and 2 body parallelograms
with rational area, and rational volume.
\begin{table}[h]
\centering
\caption{Diophantine piped with 1 rational body diagonal and rational body/face parallelograms}
{\footnotesize
\begin{tabular}{| c c c | c c c c | c c c c c c | c c c c | c c c | c c c c c c | c |}
\multicolumn{3}{c}{edges} &
\multicolumn{4}{c}{skew triangles} &
\multicolumn{6}{c}{face diagonals} &
\multicolumn{4}{c}{body diagonals} &
\multicolumn{3}{c}{face area} &
\multicolumn{6}{c}{body area} &
\multicolumn{1}{c}{vol} \\
\hline
 1 & 1 & 1 & 0 & 0 & 0 & 1 & 1 & 1 & 1 & 1 & 1 & 1 & 1 & 0 & 0 & 0 & 0 & 1 & 0 & 0 & 0 & 1 & 0 & 0 & 1 & 1 \\
\hline
\end{tabular}
}
{\small
\begin{tabular}{ r | c c c | c c c c c c | l }
 & \multicolumn{3}{c|}{Basis Vectors} & \multicolumn{6}{c}{Tetrahedron Sides} & \\
Class & $v_1$ & $v_2$ & $v_3$ & a & b & c & d & e & f & Volume\\[0.15em]
\hline\\[-0.5em]
acute & $\left(137,0,0\right)$ & $\left(\frac{3281}{137},\frac{18480}{137},0\right)$ & $\left(\frac{20475}{137},\frac{17160}{137},216\right)$ & 137 & 176 & 137 & 250 & 291 & 250 & 3991680 \Tstrut \Bstrut \\[0.5em]
\hline
\end{tabular}
}
\label{table:1bd1fa2ba}
\end{table}

\newpage

\subsubsection{Rational piped with 3 rational body diagonals, 1 rational face parallelogram and 1 rational body parallelogram}
One parallelepiped was found that had 3 rational body diagonals, 1 rational area face parallelogram and 1 rational area body
parallelogram.
\begin{table}[h]
\centering
\caption{Diophantine piped with 3 rational body diagonals and 1 rational face and 1 rational body parallelograms}
{\footnotesize
\begin{tabular}{| c c c | c c c c | c c c c c c | c c c c | c c c | c c c c c c | c |}
\multicolumn{3}{c}{edges} &
\multicolumn{4}{c}{skew triangles} &
\multicolumn{6}{c}{face diagonals} &
\multicolumn{4}{c}{body diagonals} &
\multicolumn{3}{c}{face area} &
\multicolumn{6}{c}{body area} &
\multicolumn{1}{c}{vol} \\
\hline
 1 & 1 & 1 & 0 & 1 & 0 & 0 & 1 & 1 & 1 & 1 & 1 & 1 & 1 & 0 & 1 & 1 & 0 & 0 & 1 & 0 & 0 & 0 & 1 & 0 & 0 & 0 \\
\hline
\end{tabular}
}
{\scriptsize
\begin{tabular}{ r | c c c | c c c c c c }
 & \multicolumn{3}{c|}{Basis Vectors} & \multicolumn{6}{c}{Tetrahedron Sides} \\
Class & $v_1$ & $v_2$ & $v_3$ & a & b & c & d & e & f \\[0.15em]
\hline\\[-0.5em]
acute & $\left(273,0,0\right)$ & $\left(\frac{14022}{91},\frac{20\sqrt{3552982}}{91},0\right)$ & $\left(\frac{14022}{91},\frac{3068300\sqrt{3552982}}{161660681},\frac{95760\sqrt{58624203}}{1776491}\right)$ & 273 & 431 & 442 & 431 & 442 & 560 \Tstrut \Bstrut \\[0.5em]
\hline
\end{tabular}
}
\label{table:3bd1fa1ba}
\end{table}

\subsubsection{Rational piped with 1 rational area face parallelogram and 3 rational area body parallelograms}
Computer searching found a Diophantine piped which had 1 face parallelogram with rational area, and 3 body
parallelograms with rational area.
\begin{table}[h]
\centering
\caption{Rational piped with 1 rational area face parallelogram and 3 rational area body parallelograms}
{\footnotesize
\begin{tabular}{| c c c | c c c c | c c c c c c | c c c c | c c c | c c c c c c | c |}
\multicolumn{3}{c}{edges} &
\multicolumn{4}{c}{skew triangles} &
\multicolumn{6}{c}{face diagonals} &
\multicolumn{4}{c}{body diagonals} &
\multicolumn{3}{c}{face area} &
\multicolumn{6}{c}{body area} &
\multicolumn{1}{c}{vol} \\
\hline
 1 & 1 & 1 & 0 & 1 & 0 & 1 & 1 & 1 & 1 & 1 & 1 & 1 & 0 & 0 & 0 & 0 & 0 & 0 & 1 & 1 & 0 & 1 & 0 & 1 & 0 & 0 \\
\hline
\end{tabular}
}
{\scriptsize
\begin{tabular}{ r | c c c | c c c c c c }
 & \multicolumn{3}{c|}{Basis Vectors} & \multicolumn{6}{c}{Tetrahedron Sides} \\
Class & $v_1$ & $v_2$ & $v_3$ & a & b & c & d & e & f \\[0.15em]
\hline\\[-0.5em]
obtuse & $\left(361,0,0\right)$ & $\left(-\frac{64736}{361},\frac{1320\sqrt{11041}}{361},0\right)$ & $\left(\frac{64736}{361},\frac{121178040\sqrt{11041}}{43843811},\frac{80640\sqrt{143533}}{121451}\right)$ & 361 & 663 & 424 & 425 & 424 & 448 \Tstrut \Bstrut \\[0.5em]
\hline
\end{tabular}
}
\label{table:1fa3ba}
\end{table}

\subsubsection{Rational pipeds with 3 rational face parallelograms and 2 rational body parallelograms}
Computer searching found two examples of {\it 2-ortho} {\it class} pipeds, which have all the face parallelograms rational and
2 body parallelograms rational (out of 4). These pipeds are closely related to the {\it face} {\it cuboids} found in the rational cuboid table solutions.
\begin{table}[h]
\centering
\caption{Diophantine cuboids related to {\it face} {\it cuboids}}
{\footnotesize
\begin{tabular}{| c c c | c c c c | c c c c c c | c c c c | c c c | c c c c c c | c |}
\multicolumn{3}{c}{edges} &
\multicolumn{4}{c}{skew triangles} &
\multicolumn{6}{c}{face diagonals} &
\multicolumn{4}{c}{body diagonals} &
\multicolumn{3}{c}{face area} &
\multicolumn{6}{c}{body area} &
\multicolumn{1}{c}{vol} \\
\hline
 1 & 1 & 1 & 1 & 0 & 1 & 0 & 1 & 1 & 1 & 1 & 1 & 1 & 0 & 0 & 0 & 0 & 1 & 1 & 1 & 1 & 0 & 0 & 1 & 0 & 0 & 1 \\
 1 & 1 & 1 & 0 & 1 & 1 & 0 & 1 & 1 & 1 & 1 & 1 & 1 & 0 & 0 & 0 & 0 & 1 & 1 & 1 & 0 & 0 & 1 & 0 & 0 & 1 & 1 \\
\hline
\end{tabular}
}
{\footnotesize
\begin{tabular}{ r | c c c | c c c c c c | l }
 & \multicolumn{3}{c|}{Basis Vectors} & \multicolumn{6}{c|}{Tetrahedron Sides} & \\
Class & $v_1$ & $v_2$ & $v_3$ & a & b & c & d & e & f & Volume \\[0.15em]
\hline\\[-0.5em]
2-ortho & $\left(153,0,0\right)$ & $\left(0,680,0\right)$ & $\left(0,\frac{55096}{85},\frac{17472}{85}\right)$ & 153 & 697 & 680 & 697 & 680 & 208 & 21385728 \Tstrut \Bstrut \\[0.5em]
2-ortho & $\left(185,0,0\right)$ & $\left(\frac{12593}{185},\frac{31824}{185},0\right)$ & $\left(0,0,672\right)$ & 185 & 208 & 185 & 697 & 672 & 697 & 21385728 \Tstrut \Bstrut \\[0.5em]
\hline
\end{tabular}
}
\label{table:3fa2ba}
\end{table}

\newpage

\subsubsection{Rational pipeds with 2 rational area skew triangles and 2 rational body diagonals and 1 rational area face parallelogram and 1 rational body parallelogram}
The computer searches found some interesting Diophantine pipeds which had 2 skew triangles with rational area and 2 rational body
diagonals, with 1 face parallelogram rational and 1 body parallelogram rational.

\noindent
Pipeds with rational skew area triangles were quite infrequent.
\begin{table}[h]
\centering
\caption{Diophantine pipeds with 2 rational area skew triangles}
{\footnotesize
\begin{tabular}{| c c c | c c c c | c c c c c c | c c c c | c c c | c c c c c c | c |}
\multicolumn{3}{c}{edges} &
\multicolumn{4}{c}{skew triangles} &
\multicolumn{6}{c}{face diagonals} &
\multicolumn{4}{c}{body diagonals} &
\multicolumn{3}{c}{face area} &
\multicolumn{6}{c}{body area} &
\multicolumn{1}{c}{vol} \\
\hline
 1 & 1 & 1 & 1 & 0 & 1 & 0 & 1 & 1 & 1 & 1 & 1 & 1 & 0 & 0 & 1 & 1 & 0 & 0 & 1 & 0 & 0 & 0 & 1 & 0 & 0 & 0 \\
 1 & 1 & 1 & 0 & 1 & 0 & 1 & 1 & 1 & 1 & 1 & 1 & 1 & 0 & 0 & 1 & 1 & 0 & 0 & 1 & 0 & 0 & 0 & 1 & 0 & 0 & 0 \\
 1 & 1 & 1 & 0 & 0 & 1 & 1 & 1 & 1 & 1 & 1 & 1 & 1 & 1 & 0 & 0 & 1 & 1 & 0 & 0 & 0 & 1 & 0 & 0 & 0 & 0 & 0 \\
\hline
\end{tabular}
}
{\scriptsize
\begin{tabular}{ r | c c c | c c c c c c }
 & \multicolumn{3}{c|}{Basis Vectors} & \multicolumn{6}{c}{Tetrahedron Sides} \\
Class & $v_1$ & $v_2$ & $v_3$ & a & b & c & d & e & f \\[0.15em]
\hline\\[-0.5em]
obtuse & $\left(36,0,0\right)$ & $\left(\frac{8}{3},\frac{5\sqrt{2849}}{3},0\right)$ & $\left(\frac{8}{3},-\frac{8795\sqrt{2849}}{8547},\frac{160\sqrt{1552705}}{2849}\right)$ & 36 & 95 & 89 & 95 & 89 & 160 \Tstrut \Bstrut \\[0.5em]
acute & $\left(96,0,0\right)$ & $\left(\frac{825}{8},\frac{1\sqrt{1956751}}{8},0\right)$ & $\left(\frac{825}{8},-\frac{552049\sqrt{1956751}}{15654008},\frac{2520\sqrt{16966987921}}{1956751}\right)$ & 96 & 175 & 203 & 175 & 203 & 280 \Tstrut \Bstrut \\[0.5em]
obtuse & $\left(593,0,0\right)$ & $\left(\frac{60352}{593},\frac{11040\sqrt{1533}}{593},0\right)$ & $\left(\frac{80801}{593},-\frac{788020\sqrt{1533}}{303023},\frac{620\sqrt{219219}}{511}\right)$ & 593 & 879 & 736 & 736 & 593 & 1007 \Tstrut \Bstrut \\[0.5em]
\hline
\end{tabular}
}
\label{table:2s2bd1fa1ba}
\end{table}

\subsubsection{Closest rational Diophantine (degenerate) piped to a {\it super\-perfect} piped}
The computer searches found some interesting pipeds in which all the rational checks were rational, except the volume was was 0, indicating
that these pipeds were degenerate planar figures.

\noindent
What is interesting is that these are {\it 2-ortho} {\it pipeds}.
\begin{table}[h]
\centering
\caption{Diophantine piped with all rational components except 0 volume}
{\small
\begin{tabular}{| c c c | c c c c | c c c c c c | c c c c | c c c | c c c c c c | c |}
\multicolumn{3}{c}{edges} &
\multicolumn{4}{c}{skew triangles} &
\multicolumn{6}{c}{face diagonals} &
\multicolumn{4}{c}{body diagonals} &
\multicolumn{3}{c}{face area} &
\multicolumn{6}{c}{body area} &
\multicolumn{1}{c}{vol} \\
\hline
 1 & 1 & 1 & 1 & 1 & 1 & 1 & 1 & 1 & 1 & 1 & 1 & 1 & 1 & 1 & 1 & 1 & 1 & 1 & 1 & 1 & 1 & 1 & 1 & 1 & 1 & -1 \\
\hline
\end{tabular}
}
{
\begin{tabular}{ r | c c c | c c c c c c }
 & \multicolumn{3}{c|}{Basis Vectors} & \multicolumn{6}{c}{Tetrahedron Sides} \\
Class & $v_1$ & $v_2$ & $v_3$ & a & b & c & d & e & f \\[0.15em]
\hline\\[-0.5em]
2-ortho & $\left(120,0,0\right)$ & $\left(0,182,0\right)$ & $\left(0,209,0\right)$ & 120 & 218 & 182 & 241 & 209 & 27 \Tstrut \Bstrut \\[0.5em]
\hline
\end{tabular}
}
\label{table:pluperfect}
\end{table}

\newpage

\subsubsection{Rectangular Diophantine pipeds}
There were 46 {\it rectangular} Diophantine pipeds found during the computer searches. These are the {\it body type}
of cuboid solutions. Unfortunately none of them has a rational body diagonal leading to a {\it perfect} cuboid.

\noindent
These pipeds all have the same signature for their rationality checks.
\begin{table}[h]
\centering
\caption{Rational Component checks for the {\it Rectangular} Diophantine pipeds}
{\footnotesize
\begin{tabular}{| c c c | c c c c | c c c c c c | c c c c | c c c | c c c c c c | c |}
\multicolumn{3}{c}{edges} &
\multicolumn{4}{c}{skew triangles} &
\multicolumn{6}{c}{face diagonals} &
\multicolumn{4}{c}{body diagonals} &
\multicolumn{3}{c}{face area} &
\multicolumn{6}{c}{body area} &
\multicolumn{1}{c}{vol} \\
\hline
 1 & 1 & 1 & 0 & 0 & 0 & 0 & 1 & 1 & 1 & 1 & 1 & 1 & 0 & 0 & 0 & 0 & 1 & 1 & 1 & 1 & 1 & 1 & 1 & 1 & 1 & 1 \\
\hline
\end{tabular}
}
\label{table:rectsig}
\end{table}
\begin{center}
\begin{longtable}{ r | c c c | c c c c c c }
\caption{{\it Rectangular} Diophantine Parallelepipeds}\\
\hline
Class & $v_1$ & $v_2$ & $v_3$ & a & b & c & d & e & f \\
\hline
\endfirsthead
\multicolumn{10}{c}
{\tablename\ \thetable\ -- \textit{Rectangular Pipeds - continued from previous page}} \\
\hline
Class & $v_1$ & $v_2$ & $v_3$ & a & b & c & d & e & f \\
\hline
\endhead
\hline \multicolumn{10}{r}{\textit{Continued on next page}} \\
\endfoot
\hline
\endlastfoot
rectangular & $\left(44,0,0\right)$ & $\left(0,117,0\right)$ & $\left(0,0,240\right)$ & 44 & 125 & 117 & 244 & 240 & 267 \Tstrut \Bstrut \\[0.5em]
rectangular & $\left(44,0,0\right)$ & $\left(0,240,0\right)$ & $\left(0,0,117\right)$ & 44 & 244 & 240 & 125 & 117 & 267 \Tstrut \Bstrut \\[0.5em]
rectangular & $\left(85,0,0\right)$ & $\left(0,132,0\right)$ & $\left(0,0,720\right)$ & 85 & 157 & 132 & 725 & 720 & 732 \Tstrut \Bstrut \\[0.5em]
rectangular & $\left(85,0,0\right)$ & $\left(0,720,0\right)$ & $\left(0,0,132\right)$ & 85 & 725 & 720 & 157 & 132 & 732 \Tstrut \Bstrut \\[0.5em]
rectangular & $\left(88,0,0\right)$ & $\left(0,234,0\right)$ & $\left(0,0,480\right)$ & 88 & 250 & 234 & 488 & 480 & 534 \Tstrut \Bstrut \\[0.5em]
rectangular & $\left(88,0,0\right)$ & $\left(0,480,0\right)$ & $\left(0,0,234\right)$ & 88 & 488 & 480 & 250 & 234 & 534 \Tstrut \Bstrut \\[0.5em]
rectangular & $\left(117,0,0\right)$ & $\left(0,44,0\right)$ & $\left(0,0,240\right)$ & 117 & 125 & 44 & 267 & 240 & 244 \Tstrut \Bstrut \\[0.5em]
rectangular & $\left(117,0,0\right)$ & $\left(0,240,0\right)$ & $\left(0,0,44\right)$ & 117 & 267 & 240 & 125 & 44 & 244 \Tstrut \Bstrut \\[0.5em]
rectangular & $\left(132,0,0\right)$ & $\left(0,85,0\right)$ & $\left(0,0,720\right)$ & 132 & 157 & 85 & 732 & 720 & 725 \Tstrut \Bstrut \\[0.5em]
rectangular & $\left(132,0,0\right)$ & $\left(0,351,0\right)$ & $\left(0,0,720\right)$ & 132 & 375 & 351 & 732 & 720 & 801 \Tstrut \Bstrut \\[0.5em]
rectangular & $\left(132,0,0\right)$ & $\left(0,720,0\right)$ & $\left(0,0,85\right)$ & 132 & 732 & 720 & 157 & 85 & 725 \Tstrut \Bstrut \\[0.5em]
rectangular & $\left(132,0,0\right)$ & $\left(0,720,0\right)$ & $\left(0,0,351\right)$ & 132 & 732 & 720 & 375 & 351 & 801 \Tstrut \Bstrut \\[0.5em]
rectangular & $\left(140,0,0\right)$ & $\left(0,480,0\right)$ & $\left(0,0,693\right)$ & 140 & 500 & 480 & 707 & 693 & 843 \Tstrut \Bstrut \\[0.5em]
rectangular & $\left(140,0,0\right)$ & $\left(0,693,0\right)$ & $\left(0,0,480\right)$ & 140 & 707 & 693 & 500 & 480 & 843 \Tstrut \Bstrut \\[0.5em]
rectangular & $\left(160,0,0\right)$ & $\left(0,231,0\right)$ & $\left(0,0,792\right)$ & 160 & 281 & 231 & 808 & 792 & 825 \Tstrut \Bstrut \\[0.5em]
rectangular & $\left(187,0,0\right)$ & $\left(0,1020,0\right)$ & $\left(0,0,1584\right)$ & 187 & 1037 & 1020 & 1595 & 1584 & 1884 \Tstrut \Bstrut \\[0.5em]
rectangular & $\left(234,0,0\right)$ & $\left(0,88,0\right)$ & $\left(0,0,480\right)$ & 234 & 250 & 88 & 534 & 480 & 488 \Tstrut \Bstrut \\[0.5em]
rectangular & $\left(234,0,0\right)$ & $\left(0,480,0\right)$ & $\left(0,0,88\right)$ & 234 & 534 & 480 & 250 & 88 & 488 \Tstrut \Bstrut \\[0.5em]
rectangular & $\left(240,0,0\right)$ & $\left(0,44,0\right)$ & $\left(0,0,117\right)$ & 240 & 244 & 44 & 267 & 117 & 125 \Tstrut \Bstrut \\[0.5em]
rectangular & $\left(240,0,0\right)$ & $\left(0,117,0\right)$ & $\left(0,0,44\right)$ & 240 & 267 & 117 & 244 & 44 & 125 \Tstrut \Bstrut \\[0.5em]
rectangular & $\left(240,0,0\right)$ & $\left(0,252,0\right)$ & $\left(0,0,275\right)$ & 240 & 348 & 252 & 365 & 275 & 373 \Tstrut \Bstrut \\[0.5em]
rectangular & $\left(240,0,0\right)$ & $\left(0,275,0\right)$ & $\left(0,0,252\right)$ & 240 & 365 & 275 & 348 & 252 & 373 \Tstrut \Bstrut \\[0.5em]
rectangular & $\left(252,0,0\right)$ & $\left(0,240,0\right)$ & $\left(0,0,275\right)$ & 252 & 348 & 240 & 373 & 275 & 365 \Tstrut \Bstrut \\[0.5em]
rectangular & $\left(252,0,0\right)$ & $\left(0,275,0\right)$ & $\left(0,0,240\right)$ & 252 & 373 & 275 & 348 & 240 & 365 \Tstrut \Bstrut \\[0.5em]
rectangular & $\left(275,0,0\right)$ & $\left(0,240,0\right)$ & $\left(0,0,252\right)$ & 275 & 365 & 240 & 373 & 252 & 348 \Tstrut \Bstrut \\[0.5em]
rectangular & $\left(275,0,0\right)$ & $\left(0,252,0\right)$ & $\left(0,0,240\right)$ & 275 & 373 & 252 & 365 & 240 & 348 \Tstrut \Bstrut \\[0.5em]
rectangular & $\left(351,0,0\right)$ & $\left(0,132,0\right)$ & $\left(0,0,720\right)$ & 351 & 375 & 132 & 801 & 720 & 732 \Tstrut \Bstrut \\[0.5em]
rectangular & $\left(351,0,0\right)$ & $\left(0,720,0\right)$ & $\left(0,0,132\right)$ & 351 & 801 & 720 & 375 & 132 & 732 \Tstrut \Bstrut \\[0.5em]
rectangular & $\left(480,0,0\right)$ & $\left(0,88,0\right)$ & $\left(0,0,234\right)$ & 480 & 488 & 88 & 534 & 234 & 250 \Tstrut \Bstrut \\[0.5em]
rectangular & $\left(480,0,0\right)$ & $\left(0,140,0\right)$ & $\left(0,0,693\right)$ & 480 & 500 & 140 & 843 & 693 & 707 \Tstrut \Bstrut \\[0.5em]
rectangular & $\left(480,0,0\right)$ & $\left(0,234,0\right)$ & $\left(0,0,88\right)$ & 480 & 534 & 234 & 488 & 88 & 250 \Tstrut \Bstrut \\[0.5em]
rectangular & $\left(480,0,0\right)$ & $\left(0,504,0\right)$ & $\left(0,0,550\right)$ & 480 & 696 & 504 & 730 & 550 & 746 \Tstrut \Bstrut \\[0.5em]
rectangular & $\left(480,0,0\right)$ & $\left(0,550,0\right)$ & $\left(0,0,504\right)$ & 480 & 730 & 550 & 696 & 504 & 746 \Tstrut \Bstrut \\[0.5em]
rectangular & $\left(480,0,0\right)$ & $\left(0,693,0\right)$ & $\left(0,0,140\right)$ & 480 & 843 & 693 & 500 & 140 & 707 \Tstrut \Bstrut \\[0.5em]
rectangular & $\left(504,0,0\right)$ & $\left(0,480,0\right)$ & $\left(0,0,550\right)$ & 504 & 696 & 480 & 746 & 550 & 730 \Tstrut \Bstrut \\[0.5em]
rectangular & $\left(504,0,0\right)$ & $\left(0,550,0\right)$ & $\left(0,0,480\right)$ & 504 & 746 & 550 & 696 & 480 & 730 \Tstrut \Bstrut \\[0.5em]
rectangular & $\left(550,0,0\right)$ & $\left(0,480,0\right)$ & $\left(0,0,504\right)$ & 550 & 730 & 480 & 746 & 504 & 696 \Tstrut \Bstrut \\[0.5em]
rectangular & $\left(550,0,0\right)$ & $\left(0,504,0\right)$ & $\left(0,0,480\right)$ & 550 & 746 & 504 & 730 & 480 & 696 \Tstrut \Bstrut \\[0.5em]
rectangular & $\left(693,0,0\right)$ & $\left(0,140,0\right)$ & $\left(0,0,480\right)$ & 693 & 707 & 140 & 843 & 480 & 500 \Tstrut \Bstrut \\[0.5em]
rectangular & $\left(693,0,0\right)$ & $\left(0,480,0\right)$ & $\left(0,0,140\right)$ & 693 & 843 & 480 & 707 & 140 & 500 \Tstrut \Bstrut \\[0.5em]
rectangular & $\left(720,0,0\right)$ & $\left(0,85,0\right)$ & $\left(0,0,132\right)$ & 720 & 725 & 85 & 732 & 132 & 157 \Tstrut \Bstrut \\[0.5em]
rectangular & $\left(720,0,0\right)$ & $\left(0,132,0\right)$ & $\left(0,0,85\right)$ & 720 & 732 & 132 & 725 & 85 & 157 \Tstrut \Bstrut \\[0.5em]
rectangular & $\left(720,0,0\right)$ & $\left(0,132,0\right)$ & $\left(0,0,351\right)$ & 720 & 732 & 132 & 801 & 351 & 375 \Tstrut \Bstrut \\[0.5em]
rectangular & $\left(720,0,0\right)$ & $\left(0,351,0\right)$ & $\left(0,0,132\right)$ & 720 & 801 & 351 & 732 & 132 & 375 \Tstrut \Bstrut \\
\hline
\end{longtable}
\label{table:rect}
\end{center}

\section{The Rectangular Diophantine Piped and The Perfect Cuboid}\label{sub:cuboid}

The most commonly known Diophantine parallelepiped known today is the {\it rectangular} {\it piped},
often called "the brick" or "box" or "Euler brick" also known as the {\it cuboid} or {\it integer cuboid} in math literature.
There are 7 lengths of interest in this solid, the 3 edges, the 3 face diagonals and the body diagonal.

It has been found by search, that there are 3 cases which exist, called the "body cuboid", the "face cuboid", and the "edge cuboid";
where the 7th length is irrational. This is problem D18 in {\it Unsolved} {\it Problems} {\it in} {\it Number} {\it Theory}
(UPINT) by Richard K. Guy\cite{guy}.

Walter Wyss recently attempted to show that the 7th value cannot be rational by utilizing the proof of Leonard Euler
that the sum or difference of 2 quartics cannot be a square\cite{wyss3}. He was able to parametrize the monoclinic piped and showed
that approaching a right angle {\rightangle } sets up the impossible quartic condition which Euler disproved. But he failed to prove
that no perfect integer cuboid exists, hence question is still open today.

Recently Renyxa D'Arox recently reached a milestone\cite{darox} with computer searching extending to $2^{53}$ $=$\\ $9007199254740992$ for the body diagonal, if a {\it perfect cuboid} exists, its body diagonal must exceed this value.

\section{Possible Future Studies for Diophantine Analysis}

There are some conjectures and questions created by the computer study of the 1,981,336,681 Diophantine pipeds.

Some questions:
\begin{itemize}
\item {\it Is there a rational piped with 4 rational area skew triangles?}
\item {\it Is there a perfect parallelepiped with 2 or 3 rational area face parallelograms?}
\item {\it Is there any perfect parallelepiped with 1..6 rational area body parallelograms?}
\item {\it Is there any perfect parallelepipeds with a rational area skew triangle? 2? 3? or 4?}
\item {\it Is there a perfect parallelepiped with rational volume? (asked by Sawyer, Reiter)\cite{sawyer}}
\item {\it Is there a perfect parallelepiped on the rational lattice? (asked by Sawyer, Reiter)\cite{sawyer}}
\end{itemize}

Two conjectures:
\begin{itemize}
\item {\it Conjecture 1. Other Diophantine pipeds do exist with less than 6 rational face diagonals as covered by the computer studies done here}
\item {\it Conjecture 2. No Diophantine pipeds exist with all 27 components rational}
\end{itemize}

These are just some questions and two conjectures relating to the computer search results.

Some interesting problems for further analysis are the following:
\begin{itemize}
\item {Parametrize the {\it acute} or {\it obtuse} {\it triclinic} piped as an infiinite family}
\item {Parametrize the {\it mono-orthogonal} or {\it biclinic} piped as an infinite family}
\end{itemize}

\section*{Appendix A - Proof of 5 Parallelepiped Classes}

\begin{figure}[!h]
\centering
\includegraphics[scale=1.0]{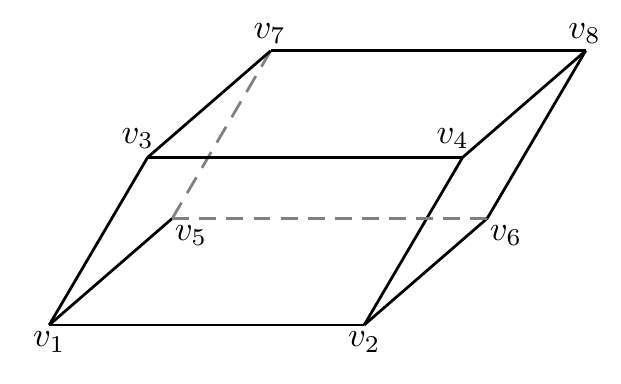}
\captionsetup{justification=centering}
\caption[caption]{Labeled vertices of the parallelepiped.}
\label{fig:vertices}
\end{figure}

Let the parallelepiped be situated so the the main vertice $v_1$ determined by the 3 basis vectors
is at the origin [0,0,0]. Figure \ref{fig:vertices} shows the parallelepiped and its 8 labeled vertices.
Previously, Figure \ref{fig:ratpiped} showed the 3 basis vectors, the same piped is shown.

The origin [0,0,0] is labeled $v_1$. The $\vec{a}$ basis vector lies along the positive
x-axis, $v_1$ to $v_2$; the $\vec{b}$ basis vector is in the x-y plane, along $v_1$ to $v_5$;
and the $\vec{c}$ basis vector is along the line from $v_1$ to $v_3$, with a non-zero z component.

We label the surface angles, as found by the angle traced by the 3 ordered vertices.

\newpage

\begin{table}[h]
\begin{center}
\caption{The 3 surface angles at each vertice of parallelepiped}
\begin{tabular}{c | c c c}
vertice & $angle_1$ & $angle_2$ & $angle_3$ \\
\hline
1 & $ v_2 v_1 v_3$ & $ v_2 v_1 v_5$ & $ v_3 v_1 v_5$ \\
2 & $ v_1 v_2 v_4$ & $ v_1 v_2 v_6$ & $ v_4 v_2 v_6$ \\
3 & $ v_1 v_3 v_4$ & $ v_1 v_3 v_7$ & $ v_4 v_3 v_7$ \\
4 & $ v_2 v_4 v_3$ & $ v_2 v_4 v_8$ & $ v_3 v_4 v_8$ \\
5 & $ v_1 v_5 v_6$ & $ v_1 v_5 v_7$ & $ v_6 v_5 v_7$ \\
6 & $ v_2 v_6 v_5$ & $ v_2 v_6 v_8$ & $ v_5 v_6 v_8$ \\
7 & $ v_3 v_7 v_5$ & $ v_3 v_7 v_8$ & $ v_5 v_7 v_8$ \\
8 & $ v_4 v_8 v_6$ & $ v_4 v_8 v_7$ & $ v_6 v_8 v_7$
\end{tabular}
\end{center}
\end{table}

We know that the interior angles on the diagonals of a parallelogram are equal,
and that 2 adjacent angles are supplementary. This fact enables us to substitute
the equivalent angles as necessary for the piped, in order to redefine all
surface angles as a subset of the 3 surface angles at vertice \#1.

\begin{center}
Equivalent angles for the parallelepiped vertices occurring on a basis plane.\vspace{6pt}\\
\begin{tabular}{c c c}
$a-b$ plane & $a-c$ plane & $b-c$ plane \\
\hline
 $v_2 v_1 v_5$ & $v_2 v_1 v_3$ & $v_3 v_1 v_5$ \\
 $v_2 v_6 v_5 = v_2 v_1 v_5$ & $v_1 v_2 v_4 = 180^\circ -v_2 v_1 v_3$ & $v_3 v_7 v_5 = v_3 v_1 v_5$ \\
 $v_1 v_5 v_6 = 180^\circ -v_2 v_1 v_5$ & $v_1 v_3 v_4 = 180^\circ -v_2 v_1 v_3$ & $v_1 v_3 v_7 = 180^\circ -v_3 v_1 v_5$ \\
 $v_1 v_2 v_6 = 180^\circ -v_2 v_1 v_5$ & $v_2 v_4 v_3 = v_2 v_1 v_3$ & $v_1 v_5 v_7 = 180^\circ -v_3 v_1 v_5$ \\
 $v_4 v_3 v_7 = v_2 v_1 v_5$ & $v_6 v_5 v_7 = v_2 v_1 v_3$ & $v_4 v_2 v_6 = v_3 v_1 v_5$ \\
 $v_4 v_8 v_7 = v_2 v_1 v_5$ & $v_5 v_6 v_8 = 180^\circ -v_2 v_1 v_3$ & $v_4 v_8 v_6 = v_3 v_1 v_5$ \\
 $v_3 v_7 v_8 = 180^\circ -v_2 v_1 v_5$ & $v_5 v_7 v_8 = 180^\circ -v_2 v_1 v_3$ & $v_2 v_4 v_8 = 180^\circ -v_3 v_1 v_5$ \\
 $v_3 v_4 v_8 = 180^\circ -v_2 v_1 v_5$ & $v_6 v_8 v_7 = v_2 v_1 v_3$ & $v_2 v_6 v_8 = 180^\circ -v_3 v_1 v_5$
\end{tabular}
\end{center}

We substitute the 3 surface angles or their supplements found at vertice \#1, for all 8 vertices,
such that all surface angles are derived from just the angles at the origin, or vertice \#1.

After substitutions, we find the surface angles.
\begin{center}
Equivalent 3 surface angles at each vertice of parallelepiped\vspace{6pt}\\
\begin{tabular}{c | c c c}
vertice & $angle_1$ & $angle_2$ & $angle_3$ \\
\hline
1 & $ v_2 v_1 v_3$ & $v_2 v_1 v_5$ & $ v_3 v_1 v_5$ \\
2 & $ 180^\circ -v_2 v_1 v_3$ & $180^\circ -v_2 v_1 v_5$ & $v_3 v_1 v_5$ \\
3 & $ 180^\circ -v_2 v_1 v_3$ & $180^\circ -v_3 v_1 v_5$ & $v_2 v_1 v_5$ \\
4 & $ v_2 v_1 v_3$ & $180^\circ -v_3 v_1 v_5$ & $180^\circ -v_2 v_1 v_5$ \\
5 & $ 180^\circ -v_2 v_1 v_5$ & $180^\circ -v_3 v_1 v_5$ & $v_2 v_1 v_3$ \\
6 & $ v_2 v_1 v_5$ & $180^\circ -v_3 v_1 v_5$ & $180^\circ -v_2 v_1 v_3$ \\
7 & $ v_3 v_1 v_5$ & $180^\circ -v_2 v_1 v_5$ & $180^\circ -v_2 v_1 v_3$ \\
8 & $ v_3 v_1 v_5$ & $v_2 v_1 v_5$ & $v_2 v_1 v_3$
\end{tabular}
\end{center}\vspace{-6pt}

In order to classify the angles, we consider the \acuteangle such that
\begin{align*}
0^\circ < \text{\acuteangle} < 90^\circ \rightarrow sign(\cos(\text{\acuteangle})) & = \phantom{-} 1 \\
	\text{\acuteangle} = 90^\circ \rightarrow sign(\cos(\text{\acuteangle})) & = \phantom{-} 0 \\
180^\circ > \text{\acuteangle} > 90^\circ \rightarrow sign(\cos(\text{\acuteangle})) & = -1
\end{align*}

We substitute these values for the surface angles in the labeling above.
For example, vertice 1 could have the signs: [-1, 0, 1], if its angles were 108\textdegree, 90\textdegree,
and 60\textdegree respectively.

We notice that if the sign of an angle = 1, then its supplement = -1, or vice versa.
If the sign of an angle = 0, then the supplement is also = 0.

Since all other 7 vertices are given as some combination of the 3 surface angles of vertice \#1,
we only have to consider the 27 permutations of -1,0,1 taken 3 at a time for vertice \#1, in order
to enumerate all possible states of the parallelepiped, and obtain an exhaustive list.

For each permutation, we create a set union of all 8 vertices to distinguish the actual unique
signature of each possible permutation group.

\begin{table}[h]
\begin{center}
\caption{The set union set of the 27 permutations for all 8 vertices of the parallelepiped}
{\scriptsize
\begin{tabular}{c|l|cc|l|cc|l}
\# & \hspace{4.0em} group & & \# & \hspace{4.0em} group & & \# & \hspace{4.0em} group \\
\hline
1 & [-1, -1, -1], [-1, 1, 1] & & 10 & [-1, -1, 0], [-1, 0, 1], [0, 1, 1] & & 19 & [-1, -1, 1], [1, 1, 1] \Tstrut \\
2 & [-1, -1, 0], [-1, 0, 1], [0, 1, 1] & & 11 & [-1, 0, 0], [0, 0, 1] & & 20 & [-1, -1, 0], [-1, 0, 1], [0, 1, 1] \\
3 & [-1, -1, 1], [1, 1, 1] & & 12 & [-1, -1, 0], [-1, 0, 1], [0, 1, 1] & & 21 & [-1, -1, -1], [-1, 1, 1] \\
4 & [-1, -1, 0], [-1, 0, 1], [0, 1, 1] & & 13 & [-1, 0, 0], [0, 0, 1] & & 22 & [-1, -1, 0], [-1, 0, 1], [0, 1, 1] \\
5 & [-1, 0, 0], [0, 0, 1] & & 14 & [0, 0, 0] & & 23 & [-1, 0, 0], [0, 0, 1] \\
6 & [-1, -1, 0], [-1, 0, 1], [0, 1, 1] & & 15 & [-1, 0, 0], [0, 0, 1] & & 24 & [-1, -1, 0], [-1, 0, 1], [0, 1, 1] \\
7 & [-1, -1, 1], [1, 1, 1] & & 16 & [-1, -1, 0], [-1, 0, 1], [0, 1, 1] & & 25 & [-1, -1, -1], [-1, 1, 1] \\
8 & [-1, -1, 0], [-1, 0, 1], [0, 1, 1] & & 17 & [-1, 0, 0], [0, 0, 1] & & 26 & [-1, -1, 0], [-1, 0, 1], [0, 1, 1] \\
9 & [-1, -1, -1], [-1, 1, 1] & & 18 & [-1, -1, 0], [-1, 0, 1], [0, 1, 1] & & 27 & [-1, -1, 1], [1, 1, 1]
\end{tabular}
\label{table:setu}
}
\end{center}
\end{table}

We discover that 5 unique classes result, after creating a set union of the 27 permutation set union groups
shown in Table \ref{table:setu}.

\begin{table}[h]
\begin{center}
\caption{Five classes of parallelepipeds}
\begin{tabular}{r|l}
name & class \\
\hline
rectangular & [0, 0, 0] \\
monoclinic & [-1, 0, 0], [0, 0, 1] \\
bi-clinic & [-1, -1, 0], [-1, 0, 1], [0, 1, 1] \\
obtuse triclinic & [-1, -1, -1], [-1, 1, 1] \\
acute triclinic & [-1, -1, 1], [1, 1, 1]
\end{tabular}
\label{table:class}
\end{center}
\end{table}
$\Box$

\par\noindent\rule{3.0in}{0.4pt}

\end{document}